\documentclass{article}
\usepackage{amssymb,amsmath,latexsym}
\usepackage{graphicx}

\begin{document}

\title{Control of a Model of DNA Division via Parametric Resonance\footnote{This work has been supported by the National Science Foundation under award number NSF-CMMI-092600.  We would like to dedicate this paper to the memory of our close collaborator Jerrold E. Marsden.
}
}

\date{September, 2012}
\author{Wang Sang Koon\footnote{Corresponding Author.  koon@cds.caltech.edu. Control and Dynamical Systems,  California Institute of Technology, Pasadena, California, 91125}, Houman Owhadi\footnote{owhadi@caltech.edu. Applied and Computational Mathematics and Control and Dynamical Systems, California Institute of Technology, Pasadena, California, 91125}, Molei Tao\footnote{mtao@cims.nyu.edu. Courant Institute of Mathematical Sciences, New York University, New York}, Tomohiro Yanao\footnote{yanao@waseda.jp. Applied Mechanics and Aerospace Engineering, Waseda University, Tokyo, Japan}}

\maketitle

\begin{abstract}
We study the internal resonance,  energy transfer,  activation mechanism, and  control of a model of DNA division via parametric resonance. While the system is robust to noise, this study shows that it is sensitive to specific fine scale modes and frequencies that could be targeted by low intensity electro-magnetic fields for triggering and controlling the division. The DNA model is a chain of pendula in a Morse potential. While the (possibly parametrically excited) system has a large number of degrees of freedom and a large number of intrinsic  time scales,  global and slow variables can be identified by (i) first reducing its dynamic to two modes
 exchanging energy between each other and (ii) averaging the dynamic of the reduced system with respect to the phase of the fastest mode. Surprisingly the global and slow dynamic of the system remains Hamiltonian (despite the parametric excitation) and the study of its associated effective potential shows how parametric excitation can  turn the unstable open state into a stable one. Numerical experiments support the accuracy of the time-averaged reduced Hamiltonian in capturing the global and slow dynamic of the full system.
\end{abstract}


{ In this paper we study the internal resonance, energy transfer,  activation mechanism,
and control of a model of DNA division via parametric resonance. While DNA macro-molecules are robust to noise, our study shows that they are sensitive to specific fine scale modes and frequencies that could be targeted by low intensity electro-magnetic fields for triggering and controlling the division.  The suggested method of control is supported not only by the observation that DNA vibrations induced by electric-fields or microwave absorbtion are  an experimental reality but also by the fact that electric field-induced molecular vibrations have already been used as a noninvasive cell transfection protocol.  Our study also raises the question on whether enzymes are using the proposed mechanism to initiate the opening of DNA strands.  This question is to put into correspondence with increasing theoretical and experimental evidence that low-frequency vibrations
 do  exist and play significant biological functions in proteins, DNA molecules,  and  other bio-macromolecules.}

\section{Introduction}
The model considered in this paper is  a chain of pendula in a Morse potential, with torsional springs between pendula, that mimic real DNA characteristics \cite{Me06, Ya04}.  Previous studies \cite{Me06, duMeMa09, EiMe07}, mainly numerical, showed that this model exhibits an intriguing phenomenon of structured activations observed in many bio-molecules:  while the system is robust to noise, it is sensitive to certain specific fine scale modes that can trigger the division.  Below, we will describe briefly the results of our analytical  study on this intriguing  phenomenon and our effort in the control of this DNA model via parametric resonance.

By using the Fourier modal coordinates \cite{duMeMa09, EiMe10, LuHoBe96}, this model can be seen as a small nonlinear perturbation of $n$ harmonic oscillators with frequencies $\omega _\alpha =\sqrt{2(1-\cos{{2\pi \alpha}/{n}})}, \alpha =0,...,n-1$.
The coarse variable of the (approximate) $0$th mode, which is the average angle of the pendula, corresponds to the angle of the frame defined by the radius of gyration in our molecular studies \cite{YaKoMaKe07, YaKoMa09, YaKoMa10}. This variable is a global and slow variable and it plays an important role as the reactive coordinate for our DNA work. Moreover, this reactive mode forms a nearly $0:1$ resonance with any other mode, each of which has an $O(1)$ frequency.  This fact leads to small denominators or coupling terms in the corresponding averaged equations or normal forms \cite{NaNa93, NaNa94, NaCh95, NaMo95, Ha99, FeLi00, TuVe03, LaZh99, BrChKiVe93}.  Since other modal frequencies are not rationally commensurate or have significant time scale separation, we do not expect strong resonance among them. Extensive simulations confirmed our expectation.  We observed: (i) the energy transfers mainly from an excited mode to the reactive mode,  triggering the division,
(ii) only an extremely small amount of  energy transfers from the excited mode to
one or two other modes via near resonance.  This observation,  together with a rigorous error estimate \cite{duMeMa09},  show that two-mode or three-mode truncation, i.e., a truncated system that includes the reactive mode, the excited mode, and perhaps a mildly affected mode should provide an adequate  reduced model for our analytical study on the activation mechanism.

By applying the method of partial averaging \cite{EiMe10, Ar98, NaNa93, NaCh95, Ha99, FeLi00, TuVe03} for nearly $0:1$ resonance, we obtained  the averaged equations for  the reduced models of this chain of Morse oscillators up to nonlinear terms of very high degree.  These averaged reduced equations not only reveal  the coupling between the action of the excited modes and the dynamics of the reactive mode, they also shed lights on the phase space structure of the activation mechanism. Moreover, they  allow us to estimate analytically the minimum activation energy for each excited mode.  These analytical estimates  not only match well with those obtained from simulations of the full DNA model, but also uncover a relationship between  the frequency of the excited mode and its corresponding  minimum activation energy.  These findings also provide an analytical and deeper understanding of the internal mechanism that is responsible for the phenomenon of structured activations.

Based on our understanding of its internal dynamics, we introduce a method for  controlling the division of this DNA model via parametric excitation, that is in resonance with its internal trigger modes.  By choosing appropriate external excitations and frictions, we  uncover a class of  trajectories that show how the parametric resonance can be used to drive the averaged reduced model from its (almost) equilibrium state to its open state.  The identification of an effective Hamiltonian (after reduction and partial averaging) opens the possibility of studying the global phase space structures of our averaged reduced model and sheds lights on the significant trajectories mentioned above.  Moreover, we extend the results for  the averaged reduced model with parametric resonance and friction  to the reduced as well as the full model with parametric resonance and friction.  These findings support the conjectures that  (i) low intensity electro-magnetic fields can be used with parametric resonance to inject energy into the trigger modes and destabilize the DNA chain and keep it near the open state for replication and transcription, (ii) enzymes may use similar method to initiate open state dynamics for DNA replication and transcription.  Although the issues of inhomogeneity, helicity, and environmental effects (such as noises)  will be ignored for now, they will be taken into consideration in our future work.


\section{A Model of DNA Division and Structured Activations}

Our analytical study has been inspired by the work of Mezi\'c and Eisenhower at UC Santa Barbara and Marsden and Du Toit at Caltech.  Their  studies \cite{Me06, duMeMa09, EiMe07, Ei09}, mainly numerical, showed that this DNA model exhibits an intriguing phenomenon of structured activations.  While Eisenhower and Mezi\'c did apply Arnold's method of partial averaging  \cite{Ar98} to a truncated Hamiltonian in their analytical study of a chain of Duffing oscillators \cite{EiMe08, Ei09, EiMe10},
they did not extend it to the DNA model.  The reason that they gave in \cite{EiMe10} for studying the Duffing case was that ``The exponential form of the Morse potential makes analytical progress difficult,...''   Since we have not seen this kind of simplification in the DNA literature, we will, in this study keep the  Morse potential (and the added difficulty).

\subsection{A Model of DNA Division}

\begin{figure}[tb]
\begin{center}
\includegraphics[width=\textwidth]{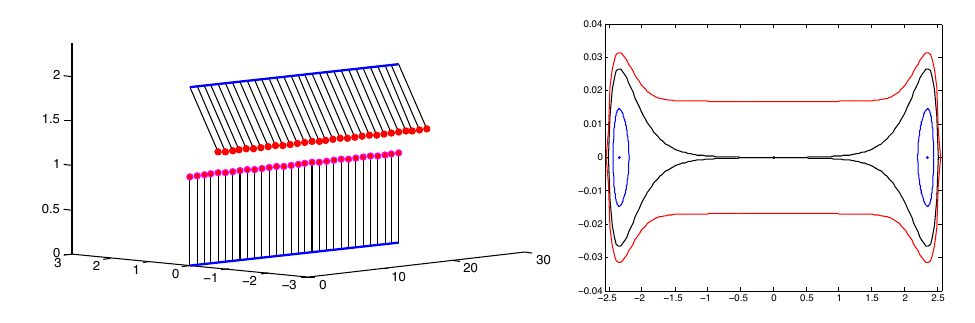}
\end{center}
\caption{\label{DNAmodel} {\footnotesize
(a)  A model of DNA division with 30 pendula.  (b)
The phase space $(\theta ,p_\theta )$ of a single pendulum in a Morse potential without coupling.
It has two stable equilibria  and a saddle.
The black curve is the homoclinic orbit  that separates  two types of motion, namely, the libration near the equilibria $(\theta _e,0)=(\pm 2.346,0)$, and the flipping across the saddle
$(0,0)$.}}
\end{figure}

The model was first introduced in \cite{Me06}.  It is a chain of equivalent pendula that rotate about the axis of a fixed backbone with an angle $\theta _k$ measured from the upward position. The pendula interact with nearest neighbors along the backbone through harmonic torsional coupling, and with pendula on the opposing immobilized strand through a Morse potential that has two stable equilibria and a saddle. See Figure \ref{DNAmodel}.

The Hamiltonian that describes the motion of these $n$ coupled pendula is given by
\begin{equation}\label{DNAmodeleq0}
H(\theta ,P_ \theta )=\sum _{k=1}^n\left[
\frac{P_{\theta k}^2}{2 m h^2}
+\frac{1}{2}S(\theta _k-\theta _{k-1})^2
+\frac{D}{2 a_d h}
\left( e^{-a_d \{ h(1+\cos \theta _k)-x_0 \} }-1\right)^2\right],
\end{equation}
with periodic boundary condition, $\theta _{0}=\theta _n$. The first term is the kinetic energy terms of $n$-pendula. The second term is the torsional coupling terms. The third term is the Morse potential terms, which model the hydrogen bonds of the respective DNA base pairs. In Eq. (\ref{DNAmodeleq0}), $m$ and $h$ represent the mass and the length of each pendulum, and $P_{\theta k} =mh^2 (d\theta_k/ dt) $ is the generalized momenta conjugate to $\theta_k$. The parameter $S$ determines the strength of the nearest neighbor coupling, while the parameter $D$ determines the strength of the Morse potential. The parameter $x_0$ determines the equilibrium distance of the Morse potential, while the parameter $a_d$ determines the width of the Morse potential. All the parameter values were chosen to best represent the typical values for the opening and closing dynamics of DNA division \cite{Chou84, Lisy96} as follows:  $m=300$ amu (typical mass of a DNA base (nucleotide)), $h=1$ nm (typical radius of DNA), $S=42$ eV, $D=0.42$ eV, $x_0 = 0.3$ nm, and $a_d = 7$ nm$^{-1}$.

After dividing the both sides of Eq. (\ref{DNAmodeleq0}) by $S$, the Hamiltonian can be non-dimensionalized as
\begin{equation}\label{DNAmodeleq}
H(\theta ,p_ \theta )=\sum _{k=1}^n\left[
\frac{1}{2}p _{\theta k}^2
+\frac{1}{2}(\theta _k-\theta _{k-1})^2
+\epsilon
\left( e^{-a(1+\cos \theta _k-d_0)}-1\right)^2\right],
\end{equation}
where $p _{\theta k} \equiv d\theta _k/ d\tau$ is the dimensionless momentum defined with respect to the dimensionless time $\tau= \sqrt{S/m h^2} t$. Thus, in the present study, one unit time ($\tau=1$) corresponds to $t = \sqrt{m h^2 /S} =$0.272 ps.  In Eq. (\ref{DNAmodeleq}), the dimensionless amplitude of Morse potential $\epsilon $ is a small parameter and is equal to  $\epsilon = D/(2 S a_d h) = 1/1400$. We have also introduced the dimensionless decaying coefficient of Morse potential $a \equiv a_d h = 7$, and the dimensionless equilibrium distance $d_0 \equiv a_d x_0/a = 0.3$.

For our analytical study, we use, as in previous  numerical studies \cite{Me06, Ei09}, a Hamiltonian system composed of 30 coupled pendula.  Figure \ref{DNAmodel}(a) is a model of 30 coupled pendula.  Before studying its dynamics, it is instructive to look at  a single pendulum in a Morse potential without coupling.  Figure \ref{DNAmodel}(b) shows the phase space of such single pendulum. It has two stable equilibria  and a saddle.  The black curve is the homoclinic orbit  that separates  two types of motion, namely, the oscillation near the equilibria $(\theta _e,0)=(\pm 2.346,0)$, and the flipping across the saddle $(0,0)$.  The n-coupled pendula have similar but much more complicated behaviors.  First, the system has two  global stable equilibria,
achieved when all the pendula have the same angular displacements $\theta _k=\theta _e$ (thus nullifying the nearest neighbors coupling) and each is positioned at the equilibria of a single pendulum.  It also has a rank one saddle at $\theta _k=0$ where all pendula are at the upward position.  For small energy, the pendula are liberating near one of the global stable equilibria  where all angles $\theta _k$ are the same and equal to $\theta _e$.  For large enough energy,  it has been observed that a local activation can cause the pendula to move collectively from one energy basin to the other and to flip across the rank one saddle.

\subsection{Phenomenon of Structured Activations}

\begin{figure}[tb]
\begin{center}
\includegraphics[width=\textwidth]{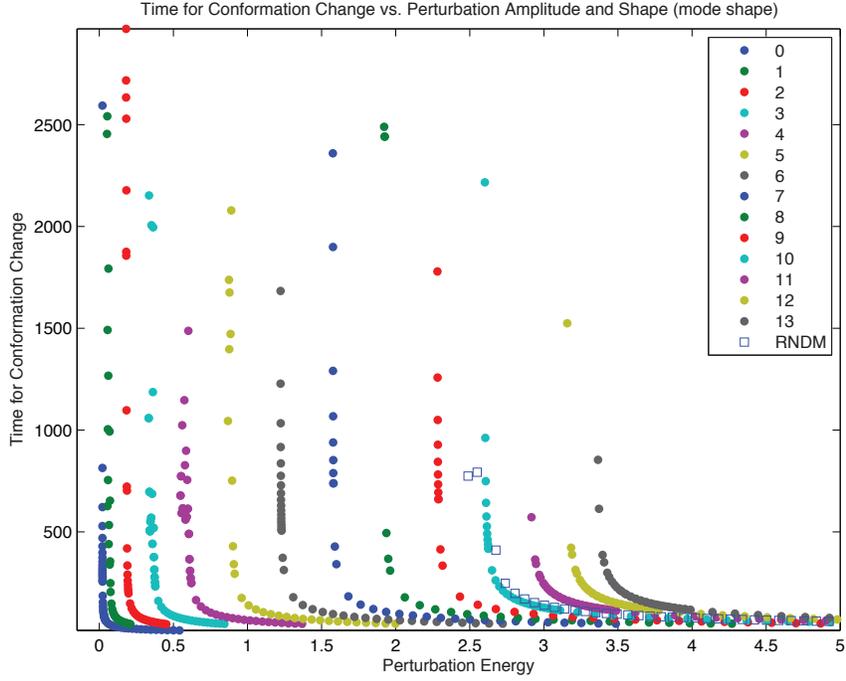}
\end{center}
\caption{\label{StructuredActivations} {\footnotesize
Figure shows the relationship between the initial amount of energy injected for various types and shapes of activation and the time for  DNA division.  Reprinted with permission from B. Eisenhower, PhD Dissertation, University of California Santa Barbara, (2009) \cite{Ei09}.
}}
\end{figure}

Previous studies  \cite{Me06, duMeMa09, EiMe07, Ei09}, mainly numerical, showed that this model exhibits an intriguing phenomenon of structured activations observed in many bio-molecules:  while the system is robust to noise, it is sensitive to certain specific fine scale modes that can trigger the division.  Figure \ref{StructuredActivations} provides the numerical data for such claim.  The figure shows the relationship between the initial amount of energy injected for various types and shapes of activation and the time for  DNA division.

In order to appreciate this figure and the phenomenon of structured activations, we need to introduce the Fourier modal coordinates $q$ which relate to the original system coordinates $\theta $ through the following linear transformation  ($\theta = Tq$):
\begin{equation}\label{Fourier}
\theta_k=\sqrt{\tfrac{2}{n}}
\sum_{\alpha =0 }^{ \tfrac{n}{2}-1}\left[
\frac{1}{\sqrt{2}}q_0
+\cos\, \frac{2\pi k\alpha }{n} q_\alpha
+\frac{(-1)^j }{\sqrt{2}} q_{\tfrac{n}{2}}
+\sin\, \frac{2\pi k\alpha }{n}q_{\tfrac{n}{2}+\alpha }
\right],
\end{equation}
where $k=1,...,n$ and $\alpha =0,...,n-1$.
Here, we have assume that the number of pendula $n$ is even.  If $n$ is odd, we merely need to have the middle column corresponding to $\alpha =\tfrac{n}{2}$ removed and the column altered accordingly.

These modal coordinates help to reveal the natural dynamics of the system
by diagonalising the linear coupling terms and rewriting the Hamiltonian as follows:
\begin{equation}
H(q,p)=\sum _{\alpha =0}^{n-1}\left(\tfrac{1}{2} p_\alpha ^2+\tfrac{1}{2}\omega _\alpha ^2 q_\alpha ^2\right)+
\epsilon \sum _{k=1}^n
U(\sum _{\beta =0}^{n-1} T_{k\beta} q_\beta )
\end{equation}
where $U(\theta )=\left( e^{-a(1+\cos \theta -d_0)}-1\right)^2$ is the Morse potential
for a single pendulum.  More specifically, the matrix $T$ is nothing but the matrix of orthonormal eigenvectors  used for the diagonalization and $\omega _\alpha $ are their corresponding  eigenvalues.  By using this coordinate system, the model can be seen as a small perturbation of $n$ harmonic oscillators  with frequencies
$\omega _\alpha =\sqrt{2-2\,\cos\, (2\pi \alpha /n)}$,  $\alpha =0,...,n-1.$
This can also be seen clearly if we write the equation of motion in the Lagrangian form:
\begin{eqnarray}\label{0to1}
\ddot q_0 \hspace{0.5in} & = & -\epsilon M_0(q_0,q_1,...,q_{n-1} ), \nonumber \\
\ddot q_\alpha + \omega _\alpha ^2 q_\alpha
&=&
-\epsilon M_\alpha (q_0,q_1,...,q_{n-1}),
\end{eqnarray}
where
\begin{equation}
M(q_0,q_1,...,q_{n-1})= \sum _{k=1}^n
U(\sum _{\beta =0}^{n-1} T_{k\beta} q_\beta)
\end{equation}
is the Morse potential term, and $M_0=\partial M/\partial q_0,
M_\alpha =\partial M/\partial q_\alpha $ are the partial derivatives of $M$ with respect to $q_0, q_\alpha $ respectively.

Observe that the  coordinate of the (approximate) $0$th Fourier mode, given as follows
\begin{equation}
q_0=\tfrac{1}{\sqrt{n}}\sum _{k=1}^n \theta _k = \sqrt {n} \bar \theta
\end{equation}
is the average amplitude of the pendula except for a constant factor of
$\sqrt{n}$.  It plays a special role as the collective variable,  the reactive coordinate, and the slow variable for the system.  Other (n-1) modal coordinates $q_\alpha $ are the fine scale  variables,  the bath coordinates, and the fast variables.

\begin{figure}[tb]
\begin{center}
\includegraphics[width=\textwidth]{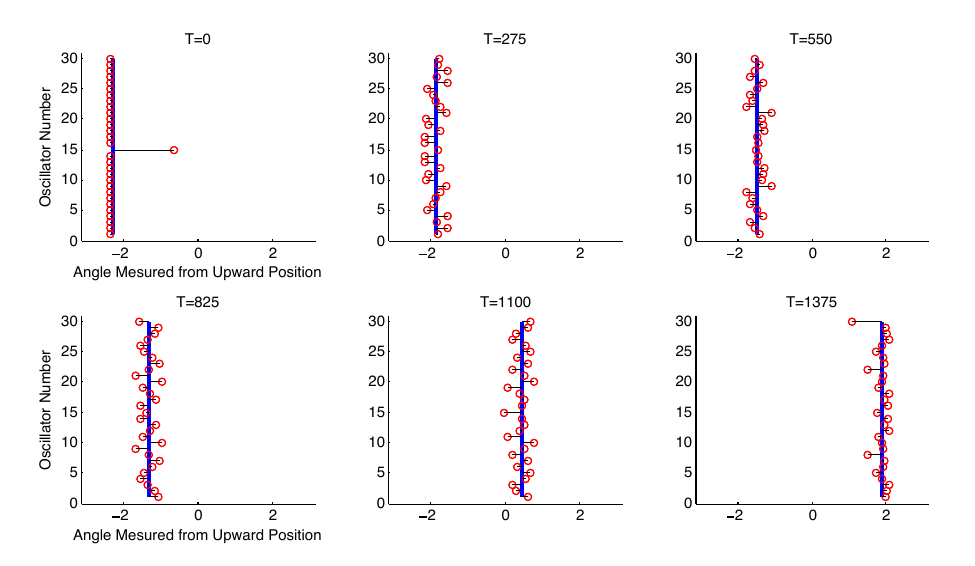}
\end{center}
\caption{\label{StructuredActivations2} {\footnotesize
Figure shows a sequence of six snapshots of the evolution of 30 pendula from one equilibrium $-\theta _e$  to the other equilibrium $\theta _e$.
Because the coupling is much stronger than the nonlinearity,
the transition is collective, closely follow the mean, the thick blue line.
}}
\end{figure}

The role of  approximate $0$th mode (or reactive) coordinate can be seen from the sequence of six snapshots of the evolution of 30 pendula from one equilibrium $-\theta _e$  to the other equilibrium $\theta _e$ of the Morse potential for a single pendulum.  For clarity of illustration, only one pendulum is perturbed as the initial activation.  Because the coupling is much stronger than the nonlinearity,
the transition is collective, closely follow the mean, the thick blue line.
Therefore, the average angle of pendula $q_0$ can be used to monitor and
mark the time of DNA division.

Now we are ready to appreciate Figure \ref{StructuredActivations}.  The figure shows the relation between the initial amount of energy injected for various Fourier modes of activation and the time for DNA division.  The curves have been  constructed by choosing the initial activation as a  single Fourier mode.  Take the magenta curve of the 4th mode as an example. Assuming that the initial position of the system is at its equilibrium  point $(q_0,...,q_{(n-1)})=(-\sqrt{n}\theta _e,0,...,0)$,  certain amount of the momentum  of the 4th mode $(p_0,...,p_{(n-1)})=(0,0,0,0,p_4(0),0,...,0)$ is injected for the initial activation at $t=0$ and its amplitude could be modified to vary  the amount of activation energy $E_{a}=p_4(0)^2/2$.  After the system evolves for awhile, the time of division $t_d$ could be determined as the time when the average angle $q_0$ first crosses the position where $q_0=0$.  In this way, a curve of $(E_a,t_d)$ is obtained that shows the amount of activation energy vs the time to DNA division for each Fourier mode.  The integration time is set for 3000 units which is approximately equal to  $2/\epsilon$.  Notice that there are 14 such curves from left to right representing those from 0th mode to 13 mode respectively.   Each curve has an asymptote at the low energy limit and will be named the minimum activation energy for each excited mode.  White ``$\square$''s show the data for  the random noise.  From the figure, we can see that the minimum activation energy depends
on the way this energy is injected into the system.  While the system is  robust to noise, it is sensitive to certain specific fine scale modes that can trigger the division.
There may also exist a relation between the minimum activation energy of each
Fourier mode and its modal frequency. In this paper, we would like to develop the analytical methods to reveal the activation mechanism and to compute the minimum activation energy for each mode. Moreover, we want to develop the techniques for controlling the real DNA division via low intensity electro-magnetic  fields and to reveal how enzymes initiate the DNA opening dynamics.


\section{Analytical Study of Structured Activations}

Let us start our analytical study of the phenomenon of structured activations.

\subsection{Nearly $0:1$ Resonance and Partial Averaging}

Recall that the equations of motion of our system are given by Eq. (\ref{0to1}) which can be seen as a nonlinear perturbation of $n$ harmonic oscillators with frequencies $\omega _\alpha =\sqrt{2-2\,\cos\, (2\pi \alpha /n)}, \alpha =0,1,...,n-1.$ Note that besides $\omega _0=0$, $\omega _\alpha $ varies from 0.2091 to 2.  Hence, the reactive mode forms a nearly  $0:1$ resonance with any other mode, each of which has an $O(1)$ frequency:
\begin{equation}
m \omega _0+0 \omega _\alpha = 0, \hspace{0.5in}
{\rm with}\hspace{0.2in} m=1.
\end{equation}
While such a relationship may not be seen as a resonance in the classical ``engineering sense'', it certainly is in the mathematical sense. This fact leads to small denominators and  coupling terms in the corresponding averaged equations or normal form. Other modal frequencies, from 1st  to 14th mode (and from 15th to 29th mode), are not rationally commensurate and do not have significant time scale separation.  We do not expect strong resonance among them.  Hence, we believe that nearly $0:1$ resonance should be the main focus of our study.

Nayfeh et al. \cite{NaNa93, NaCh95}, Haller \cite{Ha99}, Feng and Liew \cite{FeLi00}, Tuwankotta and Verhulst \cite{TuVe03}, Langford and Zhan \cite{LaZh99} and Broer et al. \cite{BrChKiVe93} have studied such a degenerate resonance.  Both Langford and Zhan and Broer et al. used the method of normal form.  Nayfeh et al., Haller, and Feng and Liew  applied a modified averaging method directly to a two mode truncation of a simple mechanical system with  parametric or external excitation.
Tukwankotta and Verhulst applied a similar method to the study of nonlinear wave
equations.

While Eisenhower and Mezi\'c did not mentioned the term,  nearly $0:1$ resonance,
in their papers,  they did apply Arnold's method of partial averaging  \cite{Ar98} to
a truncated Hamiltonian in the study of a chain of Duffing oscillators.
In \cite{EiMe10}, the reason that Eisenhower and Mezi\'c gave for studying the
Duffing case was that ``The exponential form of the Morse potential makes analytical progress difficult,...''  Since we have not seen this kind of simplification in the DNA literature, we will, in this study keep the  Morse potential (and the added difficulty).

\subsection{Two Mode Truncation Is Adequate}\label{S2Mode}

In order to carry out the analytical work for such a high dimensional system, certain reduced model is needed.  The success of Eisenhower and Mezi\'c \cite{EiMe10},
Du Toit et al. \cite{duMeMa09}, Nayfeh et al.  \cite{NaNa93, NaCh95}, and other researchers \cite{LuHoBe96, Ha99,TuVe03} shows that the method of modal truncation is a viable method if it is used with care.  Our understanding of the nearly $0:1$ resonance in our DNA model and the results of our extensive numerical simulation have convinced us that two mode truncation is adequate for our analytical study of the phenomenon of structured activations.

\begin{figure}[tb]
\begin{center}
\includegraphics[width=\textwidth]{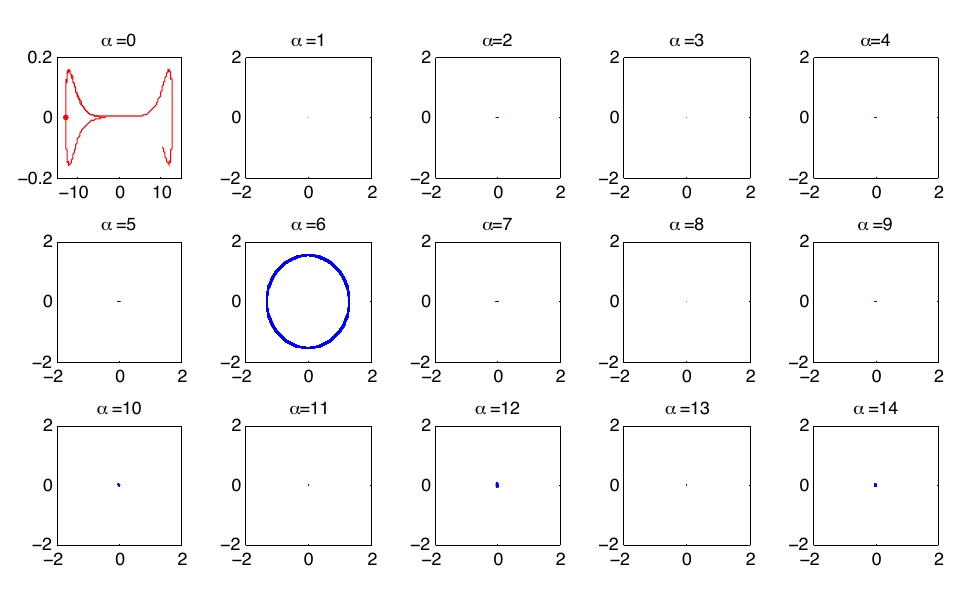}
\end{center}
\vspace{-0.2in}
\caption{\label{15Modes} {\footnotesize
The figure shows the projections of a sample trajectory on the phase spaces of the first 15 modes.  $E=1.221$ is the amount of initial kinetic energy injected into the 6th mode.
}}
\end{figure}

Figure \ref{15Modes} shows the projections of a sample trajectory on the phase spaces of the first 15 modes.  The activation has been chosen to be a single 6th mode.  From the numerical simulation which use  $q_0(0)=\sqrt{n}\theta _e$ and $p_6 (0)=\sqrt{2E}$ as the initial condition (where $E=1.221$ is the amount of energy injected into the 6th mode), we observe that (i) the energy transfer takes place over-whelmingly from the excited mode to the reactive mode, inducing the division, (ii)
only an extremely small amount of  energy transfers from the excited mode to
one or two other modes via near resonance.  Hence, we expect that the two-mode or three-mode truncations i.e., a truncated system that includes the reactive mode, the excited mode, and  perhaps a mildly affected mode, should provide an adequate reduced model for our analytical study on the nearly $0:1$ resonance, the energy transfer, and the activation mechanism.

What follow are the equations for any two-mode truncation (the reactive mode  and the $\gamma $th mode)
\begin{eqnarray}\label{EOM2Mode}
\ddot q_0 \hspace{0.5in} &=& -\epsilon M^{(2)}_0(q_0,q_\gamma )\nonumber \\
\ddot q_\gamma + \omega _\gamma ^2 q_\gamma
&=&
-\epsilon M^{(2)}_\gamma (q_0,q_\gamma)
\end{eqnarray}
where
\begin{equation}\label{ReducedMorse}
M^{(2)}(q_0,q_\gamma )= \sum _{k=1}^n
U(\sum _{\beta =\{0,\gamma \}} T_{k\beta} q_\beta)
\end{equation}
is the reduced Morse potential term, and $M_0=\partial M^{(2)}/\partial q_0,
M^{(2)}_\gamma =\partial M^{(2)}/\partial q_\gamma $ are its partial derivatives
with respect to $q_0,q_\gamma $ respectively.
Note: for notational simplicity, $M(q_0,q_\gamma )$ will be used later for $M^{(2)}(q_0,q_\gamma )$ whenever there is no ambiguity.

\subsubsection{Error Estimates for the Two Mode Truncation.} Here we show that the solution trajectories of our two mode truncation Eq. (\ref{EOM2Mode}) are within ${\mathcal O}(\epsilon)$ of the solution trajectories of the original full system described in Eq. (\ref{0to1}) for at least ${\mathcal O}(1)$ times.

The proof is an extension of the one in Du Toit et al. \cite{duMeMa09}, where a general system that has the same form as Eq. (\ref{0to1}) was studied.
In that paper, the authors (i) proposed a general technique for obtaining an $1\tfrac{1}{2}$ degree of freedom reduced system and (ii) were able to provide a rigorous error estimate for their procedure.
First, they introduced an approximation by replacing Eq. \ref{0to1}(b) with the analytical solutions of the unperturbed linear system, as defined by
\begin{equation}
Q_\alpha =A_\alpha \cos \omega _\alpha t
+\frac{B_\alpha}{\omega _\alpha} \sin \omega_\alpha  t , \hspace{0.5in} \alpha =1,...,n-1
\end{equation}
where $A_\alpha ,B_\alpha $ are the initial conditions.
Then, they obtain their reduced equation
\begin{equation}\label{duReduced}
\ddot Q_0=-\epsilon M_0(Q_0,Q_1,....Q_{n-1})
\end{equation}
via simple substitutions.  Hence, the information contained in the other modes persists in the reduced equation of the reactive mode via the initial conditions.
The error estimate that they provided claims that the solution trajectories of Eq.
(\ref{duReduced})  are within ${\mathcal O}(\epsilon)$ of the solution trajectories of the original full system described in Eq. (\ref{0to1}) for at least ${\mathcal O}(1)$ times.
This result is shown by applying a standard error analysis technique: substituting a formal expansion of the solutions, using the Lipschitz continuity of $M_0$, and then applying the Gronwall lemma as is done in the proof of Theorem 9.1 in \cite{Ver00} (see also \cite{TaOwMa11b} for a more detailed analysis on why removing small nonlinear perturbation to harmonic oscillations incurs small error).

Notice that the equations of our two mode truncation Eq. (\ref{EOM2Mode}) can be
obtained by simply setting all the initial conditions $A_\alpha ,B_\alpha $ for
Eq. (\ref{0to1}) equal to zero
except if $\alpha =\gamma$.  Hence, the solution trajectories of our two mode truncation Eq. (\ref{EOM2Mode}) should  also be within ${\mathcal O}(\epsilon)$ of the solution trajectories of the original full system described in Eq. (\ref{0to1}) for times ${\mathcal O}(1)$.

Moreover, this error estimate is valid for an arbitrary forcing function.  But for our DNA model, the exponential decay of the Morse potential with the distance implies that when the pendula escapes the immediate vicinity of the opposing pendula, the Morse potential and its consequent perturbation are effectively zero.  Hence, our error bound is very loose because it utilizes only the fact that $\epsilon$ is small, but not the specific feature of our DNA model, which is that the forcing term is also small in a large region of phase space.  This is why we numerically observed that the two mode truncation remains accurate over a timescale much longer than $\mathcal{O}(1)$.

\begin{figure}[tb]
\begin{center}
\includegraphics[width=\textwidth]{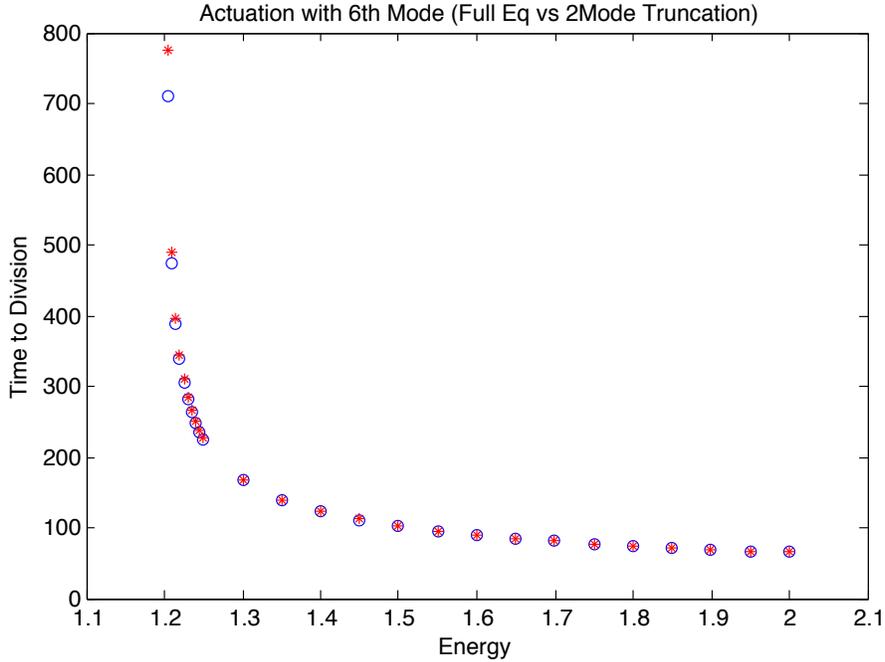}
\end{center}
\caption{\label{MAE} {\footnotesize
Figure shows the amount of activation energy vs the time of DNA division for the 6th mode.
The blue ``o''s are the data points for the full equations of 30 modes;
the red ``*''s are correspondent data points for the two mode truncation.
Both have almost the same  minimum activation energy $1.205$.
}}
\end{figure}

Figure \ref{MAE} show the amount of activation energy vs the time of DNA division for the 6th mode.
The blue ``o''s are the data points for the full equations of 30 modes;
the red ``*''s are correspondent data points for a two mode truncation (the reactive mode and the 6th mode).
Both have almost the same  minimum activation energy $1.205$.
Even their time of division differs very little if the actuation energy is slightly larger than the minimum activation energy (away from the asymptote).
The above observation also holds true for all the other modes.

Now we are ready to do the truncation and apply the method of partial averaging to obtain the
averaged equations for  the reduced models
of this chain of coupled Morse
oscillators, and use them to study the activation mechanism and
compute the minimum activation energy for each mode.

\subsection{Averaged Reduced Hamiltonian}

Since in our case, the partial averaging of Euler-Lagrangian equations is
equivalent to the partial averaging of its Hamiltonian and our main concern is on energy transfer, we preferred to
use the Hamiltonian formulation in this study.
The equivalency was shown in \S \ref{Equivalent} of the Appendix.

Recall that  the Hamiltonian
of the two mode reduced model is given by
\begin{equation}
H_2^\gamma (q,p)=\sum _{\alpha =\{0,\gamma \}}\left(\tfrac{1}{2} p_\alpha ^2+\tfrac{1}{2}\omega _\alpha ^2 q_\alpha ^2\right)+\epsilon \sum _{k=1}^n
U(\sum _{\beta =\{0,\gamma \}}T_{k\beta} q_\beta ).
\end{equation}
The averaged reduced Hamiltonian can be obtained as follows.
First, we approximate (via Taylor expansion) the Morse potential $U$, which involves the exponential function, by a polynomial of degree 26 at $\theta =0$ \begin{equation}
H_2^\gamma (q,p)=\sum _{\alpha =\{0,\gamma \}}
\left(
\tfrac{1}{2} p_\alpha ^2+
\tfrac{1}{2}\omega _\alpha ^2 q_\alpha ^2\right)
+\epsilon \sum _{k=1}^n
\sum _{j=0}^{26}a_j(\sum _{\beta =\{0,\gamma \}}T_{k\beta} q_\beta )^j
\end{equation}
where $a_j, j=0,..., 26$ are the coefficients of the expansion.
As pointed out in
\cite{EiMe10}, the Morse potential is indeed a difficult function to work with.  Since
the geometry of the Morse potential requires a polynomial of degree at least $26$ for an accurate approximation,
we decide to use such a high degree expansion for our computation of the averaged reduced Hamiltonian.

Then, we use the angle-action coordinates defined as follows
\begin{equation}
q_\gamma =\sqrt{2I_\gamma /\omega _\gamma } \sin \phi _\gamma,
\hspace{0.5in} p_\gamma =\sqrt{2I_\gamma \omega _\gamma } \cos \phi _\gamma .
\end{equation}
and rewrite the reduced Hamiltonian as
$H_2^\gamma (q_0,p_0,I_\gamma ,\phi _\gamma)$.
Notice that, besides $q_0,p_0$, the action $I_\gamma $ is also a slow variable.
Hence, the averaged reduced Hamiltonian can be obtained by
averaging the only fast variable $\phi _\gamma$:
\begin{equation}
\bar H_2^\gamma  (\bar q_0,\bar p_0, \bar I_\gamma)
=\tfrac{1}{2\pi}\int _0^{2\pi} H_2^\gamma (q_0,p_0,I_\gamma ,\phi _\gamma)d\phi _\gamma .
\end{equation}

After renaming variables, the averaged reduced Hamiltonian is given by
\begin{equation}\label{AverageHam}
\bar H_2 =\frac{1}{2}y^2+\omega I+\epsilon \left(
na_0+\sum _{k=0}^{13} c_{2k}(I)x^{2k}\right)
\end{equation}
where
$x=\bar q_0,y=\bar p_0$ are the Cartesian coordinates of the reactive mode;
$I, \omega $ are the action and the constant frequency of the other mode;
$\epsilon na_0=0.0214$ is the energy value at  the saddle;
$c_{2k}(I)$ are polynomials in $I$.  For example, $c_0(I)$ is a 13 degree polynomial in $I$ given by
\begin{equation}
c_0(I)=\sum _{j=1}^{13} b_j I^j
\end{equation}
where $b_j$ are numerical coefficients.  Also, for the simplicity of notation, the letter $\gamma $ is dropped from the averaged reduced Hamiltonian $\bar H_2^\gamma $.

\subsubsection{Error Estimates for the Averaging.}  By applying the standard averaging theory for differential equations (see for instance \cite{SaVeMu10}; we again recall that averaging the Hamiltonian is equivalent to averaging the equations), we can conclude that the solution trajectories of the averaged reduced equations are within ${\mathcal O}(\sqrt{\epsilon})$ of the solution trajectories of the reduced equations for at least ${\mathcal O}(1/\sqrt{\epsilon})$ times. 

\subsection{Error Estimates of the Entire Treatment.}
There are three approximations that we made in our treatment: Taylor approximations of the potential, truncation into two modes, and averaging. We found that the truncation and averaging respectively induces $\mathcal{O}(\epsilon)$ and $\mathcal{O}(\sqrt{\epsilon})$ errors for at least $\mathcal{O}(1)$ and ${\mathcal O}(1/\sqrt{\epsilon})$ times. In addition, as long as the solution remains bounded, Taylor approximations result in a small $o(\epsilon)$ error in the nonlinear forces, which again by Gronwall only induces $\mathcal{O}(\epsilon)$ error in the solution till at least $\mathcal{O}(1)$ times. Put together, our entire treatment induces at most an $\mathcal{O}(\sqrt{\epsilon})$ error till at least $\mathcal{O}(1)$ time.

Moreover, this error estimate is valid for an arbitrary forcing function.  But for our DNA model, the exponential decay of the Morse potential with the distance implies that when the pendula escapes the immediate vicinity of the opposing pendula, the Morse potential and its consequent perturbation are effectively zero.  Hence, our error bound is very loose because it utilizes only the fact that $\epsilon$ is small, but not the specific feature of our DNA model, which is that the forcing term is also small in a large region of phase space.  This is why we numerically observed that our approximation remains accurate over a timescale much longer than $\mathcal{O}(1)$.

\subsection{Phase Space of Averaged Reduced Equations}

Given the averaged reduced Hamiltonian, we obtain the averaged reduced Hamiltonian equations as follow
\begin{eqnarray}\label{AverageReduced}
\dot x&=&y,
\hspace{2.2in}   \dot I=0 \nonumber \\
\dot y&=&-\epsilon \left(\sum _{k=1}^{13}2kc_{2k}(I)x^{2k-2}\right)x;
\hspace{0.45in}
\dot \phi =\omega + \epsilon\left(\sum_{k=0}^{13}\frac{dc_{2k}}{dI}x^{2k}\right).\end{eqnarray}
Notice that the action $I$ of the excited mode is a constant of motion.  Hence, the averaged equations
of the reactive mode, i.e., the $x,y$ equations, and
their phase space structures are parametrized by $I$.

\begin{figure}[tb]
\begin{center}
\includegraphics[width=\textwidth]{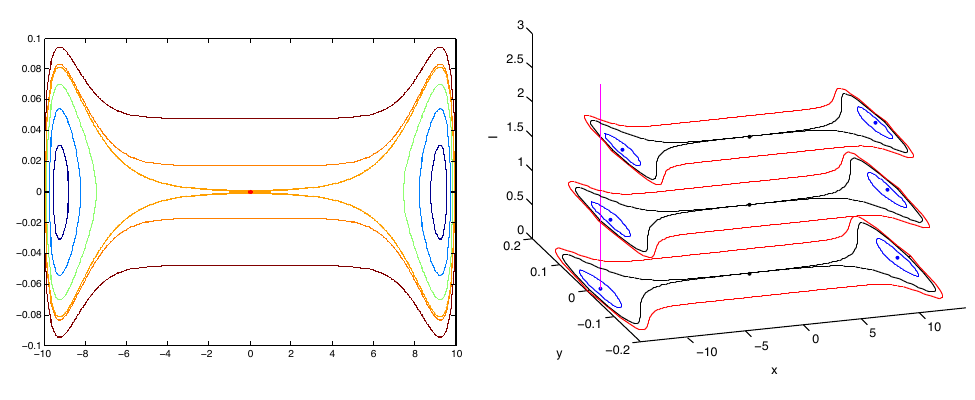}
\end{center}
\caption{\label{PhaseSpace} {\footnotesize
(a) The contour plots for the averaged reduced Hamiltonian when  $I=5$.
This phase space $(x, y)$  has a separatrix which separates two types of motion: liberation (also referred to as ``breathing'' in \cite{Chou84}) and flipping.
(b) Three such phase spaces stacked up in increasing $I$ in the
$(x,y,I)$ space.
Notice that the phase space and the separatrix ``shrinks'' towards the saddle as
$I$ increases. Together, these separatrices  form a homoclinic manifold and can be used in studying the minimum activation energy for each excited mode.
}}
\end{figure}

Therefore, the averaged reduced Hamiltonian, Eq. (\ref{AverageHam}),
can be used to study the phase spaces of the reactive mode of the averaged reduced equations that are parametrized by $I$.
Figure \ref{PhaseSpace}(a) shows the  contour plots of the averaged reduced Hamiltonian and  the phase space for the reactive mode of the
averaged reduced equations, when  $I=5$.
Notice that this phase space has a
separatrix which separates two types of motion: liberation and flipping.
Figure \ref{PhaseSpace}(b) shows three of such
phase spaces stacked up in increasing $I$ (action of excited mode).
In the $(x,y,I)$ space, the vertical axis can also be seen as
the axis of increasing activation
energy by scaling with $\omega$, $E_{act}=\omega I$.
Notice that the phase space and the separatrix ``shrinks'' towards the saddle as $E_{act}$ increases. Together, these separatrices  form a homoclinic manifold and can be used in studying the minimum activation energy for each excited mode.

\subsection{Analytical Study of Minimum Actuation Energy}\label{Analytical}

Recall that for large enough energy in the excited mode, the energy transfered  to
the reactive mode will
surpass those at the saddle and the system will be driven from
its initial equilibrium across the separatrix, causing the flipping.
See Figure \ref{15Modes}.
But for the averaged reduced system, this process manifests itself via the changes in the phase space and the separatrix of the reactive mode parametrized by $I$.
See Figure \ref{PhaseSpace}(b) for illustration.
For activation energy slightly larger than the minimum activation energy,
$E_{min} (=\omega I_{m}$),
the initial stable equilibrium at $(x_{e},0)$, now parametrized
by $I_{min}^+$,
will cross the  separatrix  parametrized by $I_{\rm m}^+$,
move into  the flipping region of phase space parametrized by $I_{m}^+$ and induce the DNA division.
Hence, $E_{min}$ can be found by the condition that
the point $(x_{e},0,I_{m})$ is on the separatrix that passes through $(0,0,I_{m})$.
Since the separatrix is the energy curve $(x,y)$ defined by
\begin{equation}
\bar H_2(x,y,I_{m})=\bar H_2(0,0,I_{m})
\end{equation}
where $\bar H_2$ is the averaged reduced Hamiltonian, $E_{min}$ can be found by solving the following equation for $I_{m}$
\begin{equation}
\bar H_2(x_{e},0,I_{m})=\bar H_2(0,0,I_{m})
\end{equation}
where $(0,0,I_m)$ is a saddle at $I_m$ and $(x_e,0)$ is the initial equilibrium point
at $I=0$.

As pointed out earlier, $\bar H_2$ of Eq. (\ref{AverageHam}) has been computed using a Taylor expansion of the Morse potential $U$ at $\theta =0$.
Strictly speaking, it may be better to denote it as $\bar H_2^0$.
Even though we have already tried to find an excellent Taylor expansion that can take care of the approximation of $U$ at both the neighborhood near the saddle
$\theta =0$ and the basins around the stable equilibria $\theta =\pm \theta _e$,  the accuracy of $U$ and $H_2^0$ at the latter is still not as good as those at the former.
For example,
while the actual stable equilibria should be at
$(\pm \sqrt{n}\theta _e,0)=(\pm 12.85,0)$,
the approximated ones will be at
$(\pm 12.59,0)$ if $\bar H_2^0$ is used.
Since any significant error in the approximation of the Morse potential $U$
at the stable equilibria $\pm \theta _e$
will cause a large error in the estimation of the minimum activation energy, another averaged Hamiltonian $H_2^e$ has also been computed
using an expansion of the Morse potential $U$ at one of the stable equilibria
$\theta =-\theta _e$:
\begin{equation}
\bar H_2^e=\frac{1}{2}y^2+\omega I+\epsilon
\left(\sum _{j=0}^{26} c'_{j}(I)(x+\sqrt{n}\theta _e)^{j}\right).
\end{equation}
Here, the letter $e$ stands for the word ``equilibria''.  Now,
The minimum activation energy can be analytically estimated by solving
the following equation
\begin{equation}
\bar H_2^e(-\sqrt{n}\theta _e,0,I_m)=\bar H_2^0(0,0,I_m).
\end{equation}
Straightforward substitution will give us
\begin{equation}
\omega I_m+\epsilon  c'_{0}(I_m)
=\omega I_m +\epsilon (na_0+c_0(I_m)).
\end{equation}
After simplification, the minimum activation energy can be analytically estimated by
solving
\begin{equation}
c_0'(I_m)-c_0(I_m)-na_0=0.
\end{equation}
For the two-mode truncation of the reactive mode and the 6th mode, $I_m/\omega _6=0.8539$.
Therefore, $E_{min}=\omega _6\times 0.8539\times \omega _6=1.1801$, as compared with $1.205$ from numerical simulation.

\begin{figure}[tb]
\begin{center}
\includegraphics[width=\textwidth]{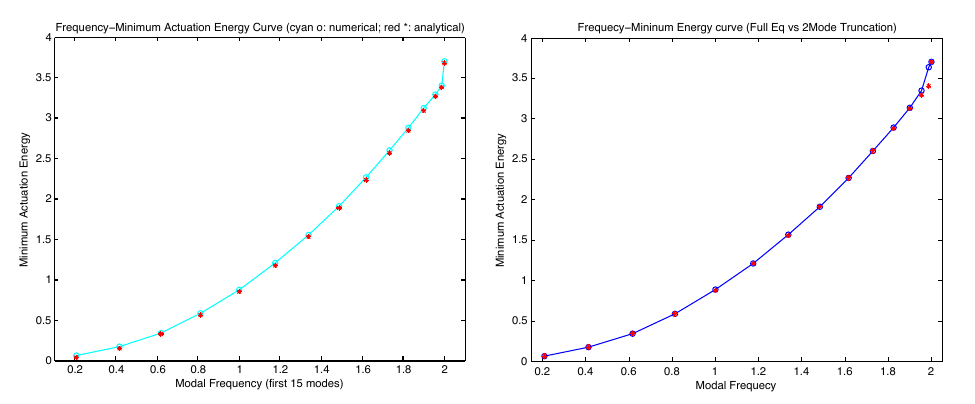}
\end{center}
\caption{\label{MAECompare} {\footnotesize
(a) Analytical estimation of minimum activation energies vs numerical simulation obtained from two mode reduced models.
The data for  magenta  ``*''s are from analytical computation.
The data for cyan ``o''s
are from numerical simulation obtained from the two-mode truncations.
They match very well.
(b) Numerical simulation of two mode reduced models vs numerical simulation obtained from the full model.
The blue ``o''s are the data from numerical simulation of the full equations (30 modes).
The red ``*'' are from simulation of the two-mode
truncations.
}}
\end{figure}

Figure \ref{MAECompare}(a) compares the analytical estimation of minimum activation energy with the numerical simulation obtained from two-mode reduced model.
The data for cyan ``o''s
are from the numerical simulation obtained from the two-mode truncation.
For example, the 6th data point (for the 6th mode) with
$\omega _6 =1.17$ needs minimum activation energy
${\rm E}^{(6)}_{m}=1.205$.
The data for  magenta  ``*''s are from analytical computation.
For example, the 6th data point with
$\omega _6=1.17$ has $E^{(6)}_{min}=1.1801$
Clearly, the data from numerical simulation and the data from analytical computation match
very well.
Moreover, $E_{min}$ (except the 15th mode) can be approximated as a parabola
\begin{equation}
E_{min}=0.8539\times \omega ^2.
\end{equation}
Here, we would like to remark that the use of two different Taylor expansions of the Morse potential function $U$, one at $\theta =0$ and another at $\theta =\theta _e$, not only allows us to deal with the difficulty posed by the exponential form of the Morse potential function but also improves the accuracy of our analytical estimation
of the minimum activation energy for each excited mode.

Figure \ref{MAECompare}(b) shows the comparison between numerical simulation obtained from the two mode truncations with the numerical simulation obtained from the full model.   The blue ``o''s are the data from numerical simulation of the full equations (30 modes).
For example, the 6th data point (for the
6th mode) with $\omega _6= 1.17$ needs minimum activation enesrgy
$ {\mathcal E}^{(6)}_{m}=1.205$.
The red ``*'' are from simulations of the two-mode
truncations ${\rm E}^{(6)}_m=1.205$.
Besides the 14th mode (with error less than $5\%$),
all other data values of blue ``o''s and red ``*''s have differences
less than 1\%.
Therefore, for the study of the minimum actuation energy, the analytical estimation  provides accurate prediction for the full system.
.

\subsection{Summary of Analytical Results on Structured Activations}

By applying the method of partial averaging, we have obtained the
averaged reduced equations
Eq. (\ref{AverageReduced}) for a chain of coupled Morse oscillators.
(i) These equations reveal the coupling between  the action and energy of the excited mode and the  dynamics of the reactive mode, as well as the phase space structure of the activation mechanism.
(ii) They allow  us to estimate analytically the minimum activation energy for each excited mode
and discover a relation between the  frequency of the excited mode and their corresponding  minimum activation energy.
(iii)  These estimates
match very well with the numerical
simulations obtained from the reduced and  the full model.
The results show that the  nearly $0:1$ internal resonance is responsible for the structured activation of our DNA model.

\subsection{Remark on  Pitchfork Bifurcation}

We observe that the reactive mode of the averaged reduced model has a pitchfork bifurcation at $I_b \approx 38$.
This can be seen clearly by studying the following averaged equations of the reactive mode parametrized by $I$:
\begin{equation}\label{AverageReduced2}
\dot x=y,
\hspace{1in}
\dot y=-\epsilon \left(\sum _{k=1}^{13}2kc_{2k}(I)x^{2k-2}\right)x;
\end{equation}

\begin{figure}[tb]
\begin{center}
\includegraphics[width=\textwidth]{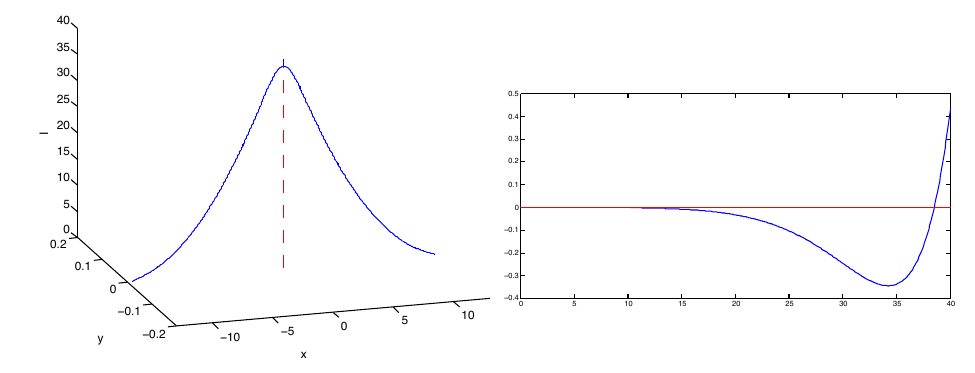}
\caption{\label{Pitchfork} {\footnotesize
Figure (a) shows that the reactive mode of the averaged reduced system has a pitchfork bifurcation at the parameter value $I\approx 38$.
Figure (b) shows a graph of $c_2(I)$ which has a simple zero at $I\approx 38$.
}}
\end{center}
\end{figure}

Notice that the equilibria of this system are given by $(0,0)$ and $(x_e(I),0)$ where
$x_e(I)$ are the solutions of
\begin{equation}
\sum _{k=1}^{13}2kc_{2k}(I)x^{2k-2}=0.
\end{equation}
See the blue curve
in Figure \ref{Pitchfork}(a).
Clearly, $(x_e(I),0)$ are stable equilibria, except at $(0,0)$.  In fact, they are surrounded by liberation
contours as shown in Figure \ref{PhaseSpace}.  As for the equilibria $(0,0)$ for each $I$,  their stability are determined by the sign of the coefficient of the linear term in
$x$, namely, $-2\epsilon c_2(I)$.  For $c_2(I)<0$, they are saddles.
For $c_2(I)>0$, they are stable equilibria.  Figure \ref{Pitchfork}(b) shows a graph of
$c_2(I)$ which has a simple zero at $I_b \approx 38$.  Hence, the system has a pitchfork
bifurcation as shown in Figure \ref{Pitchfork}(a).  As will be shown in the following section,
this bifurcation will play an important role in our study on the control of this DNA model via parametric resonance.


\section{Control via Parametric Resonance}

Building on our understanding of its internal dynamics, we want to
control the division of this DNA model via parametric excitation,
that is in resonance with its internal trigger modes.
This effort is guided by two observations and two conjectures:
The two observations are (i)
the averaged reduced model has a pitchfork bifurcation at
$I_b \approx 38$ ---- the physical
interpretation of this is that if enough energy is injected into the trigger mode
with a value larger than its bifurcation value, the chain of pendula will remain near its open state;
(ii) parametric resonance is an efficient way for energy transfer from an external source.  Hence,
our two conjectures are  (i)
low intensity electro-magnetic  fields can be used with the mechanism of parametric
resonance to pump energy into the trigger modes and keep the real DNA chain near the opening state for replication and transcription.  Here we note that electric field-induced molecular vibrations have already led to the development of a noninvasive
cell transfection protocol that enables foreign DNA
molecules   to cross cell membranes and penetrate into
the cytoplasm by eliciting vigorous vibration between
molecules and cells \cite{Zhang89, Song2004};
(ii) enzymes may use similar method to initiate the opening  dynamics
for DNA replication and transcription. For discussions and evidence of important relations between low-frequency vibrations in  DNA molecules (and other bio-macromolecules) and significant biological functions, we refer to  \cite{Chou84}, \cite{Chou84b}, \cite{Careri75}, \cite{Sobell79}, \cite{Sobell83}, \cite{Englander80}, \cite{Zhou81}, \cite{Letellier1986}.

Our three main results in this section are as follows:
(i) We identify excitation parameters  (as a function of friction) and a class of  trajectories that show how  parametric resonance can be used to  drive the averaged reduced model from its almost equilibrium state to its open state.
(ii) By identifying the averaged reduced effective Hamiltonian (after parametric excitation), we uncover the global phase space structure
and analyze the characteristic trajectories mentioned above.
(iii) We extend the results for the  averaged reduced model with parametric excitation and friction to the reduced and to the full model with parametric resonance and friction. The findings uncover a method  for  controlling  DNA division via
parametric resonance.
Although the issues of inhomogeneity, helicity, and environmental effects (such as noises) are ignored for now, they will be taken into consideration in our future work.

\subsection{Equations of Motion with Parametric Excitation}

The equations of motion of our full DNA model with parametric excitation and frictions can be written as follows:
\begin{equation}\label{Full_PR}
\ddot \theta _k -  (\theta _{k+1} -2\theta _k+\theta _{k-1}) -\epsilon U'(\theta _k)=
\epsilon \theta _k f\cos {\Omega t}
-\epsilon \mu \dot \theta _k
\end{equation}
where $k=1,...,n$, $\theta _{0}=\theta _n$, $U$ is the Morse potential function and $U'$ is its derivative.
Notice that the left hand side is the original equations of motions, without parametric excitation or frictions.  See Eq. \ref{DNAmodeleq}.  Moreover,
$f,\Omega$ are the amplitude and the frequency of the parametric excitation, respectively;
$\Omega $ is in a nearly 1:2 parametric resonance with a chosen internal trigger mode $\omega _\gamma $;
$\mu $ is the frictional coefficient.

Note the friction terms represents energy loss caused by interaction with surrounding molecules. Noise terms (ignored here) would represent energy gain caused by (thermal) interaction with surrounding molecules \cite{Muto90}.
The parametric excitation terms represent interactions with electro-magnetic fields, we note that DNA vibrations induced by electric-fields or microwave absorbtion are an experimental reality \cite{Zhang89, Song2004}.

As pointed out in Ref. \cite{Zhang89}, the bases of DNA have dipole moments \cite{Devoe62}, which could couple with an external electromagnetic field. In principle, this coupling can induce a periodic wave-like forcing over DNA bases. Since an electromagnetic field has a polarization, it is reasonable to assume that the coupling between an electromagnetic field and a DNA base is dependent on the orientation of the base, $ \theta _k$. This motivates us to introduce the parametric excitation term in our model equation (the first term on the right-hand side of Eq. (\ref{Full_PR}) above).

After using the Fourier modal coordinates, the equations  of motion for a corresponding two mode reduced model are given by
\begin{eqnarray}\label{Reduced_PR}
\dot q_0&=&p_0
\hspace{1in}
\dot p_0=-\epsilon M_0 -\epsilon  \mu p_0+\epsilon  fq_0\cos \Omega t. \nonumber\\
\dot q_\gamma &=&p_\gamma
\hspace{1in}
\dot p_\gamma =-\omega _\gamma ^2q_\gamma -\epsilon M_\gamma
-\epsilon  \mu p_\gamma +\epsilon  fq_\gamma \cos \Omega t
\end{eqnarray}
where $q_0,p_0 $ and $q_\gamma ,p_\gamma $ are the coordinates for the reactive mode and the $\gamma $-mode, respectively; $M$ is the reduced Morse potential term and $M_0, M_\gamma $ are its partial derivatives with respect to $q_0, q_\gamma $ respectively.  See Eq. (\ref{EOM2Mode}) and Eq. (\ref{ReducedMorse}) for detail.

After applying the method of partial averaging to Eq. (\ref{Reduced_PR}) using the
angle-action variables
\begin{equation}\label{AngleAction2}
q_\gamma =\sqrt{4I/\Omega }\sin \left(\tfrac{\Omega }{2}t+\beta \right),
\hspace{0.5in}
p_\gamma =\sqrt{I\Omega }\cos \left(\tfrac{\Omega }{2}t+\beta \right),
\end{equation}
renaming variables, and setting $\Omega =2\omega$,
we obtain the averaged reduced equation of motion
\begin{eqnarray}\label{ARE_PR}
\dot x&=&y,
\hspace{1.35in}
\dot \beta =\epsilon (\bar M_I
-\sigma/2\omega
+f\cos 2\beta/4\omega),\nonumber \\
\dot y&=&-\epsilon (\bar M_x
+ \mu y),
\hspace{0.5in}
\dot I=
\epsilon I( f \sin 2\beta/2\omega-\mu),
\end{eqnarray}
where $I,\beta $ are the action and the phase of the chosen internal trigger mode; $\sigma $ is a detuning parameter, defined by
$\Omega ^2/4=\omega^2+\epsilon \sigma $; $\bar M$ is the averaged Morse term, and
$\bar M_x, \bar M_I$ are its partial derivatives with respect to $x,I$ respectively.
See  Eq. (\ref{AverageReduced}) for comparison.

The derivation is similar to what has been done in \cite{NaCh95, Ha99} and is
parallel to the method used in  \S \ref{Equivalent} of the Appendix
with only minor modification.  We first rewrite Eq. (\ref{Reduced_PR})
in the following first order form:
\begin{eqnarray}\label{Reduced_PR2}
\dot q_0&=&\hspace{0.13in}\sqrt{\epsilon}p_0 \hspace{1.05in} q_\gamma =p_\gamma
\nonumber \\
\dot p_0&=&-\sqrt{\epsilon}R_0 \hspace{1in} p_\gamma =
-\omega ^2_\gamma q_\gamma -\epsilon R_\gamma .
\end{eqnarray}
where
\begin{eqnarray}\label{R0}
R_0&=& M_0+\mu p_0-fq_0 \cos \Omega t \nonumber\\
R_\gamma &=& M_\gamma +\mu p_\gamma -fq_\gamma \cos \Omega t .
\end{eqnarray}
By using the
angle-action variables defined by Eq. (\ref{AngleAction2}),
we can transform
Eq. (\ref{Reduced_PR2}) into the Lagrangian Standard Form
\begin{eqnarray}\label{Reduced_PR3}
\dot q_0&=&\hspace{0.13in}\sqrt{\epsilon}p_0 \hspace{1.05in}
\dot \beta  =\hspace{0.13in}
\epsilon (\tilde R_\gamma  /\sqrt{I \Omega })\sin \phi
\nonumber \\
\dot p_0&=&-\sqrt{\epsilon}R_0 \hspace{1in}
\dot I =-\epsilon \tilde R_\gamma \sqrt{4I /\Omega  }\cos \phi
\end{eqnarray}
where $\phi = (\tfrac{\Omega }{2} t +\beta )$ and
\begin{equation}
 \tilde R_\gamma =R_\gamma -\sigma \sqrt{4I/\Omega } \sin \phi .
 \end{equation}
Hence,  we can apply the standard  averaging theory (by averaging $t$ from $0$ to $4\pi/\Omega  $) and obtain the averaged reduced equations
\begin{eqnarray}\label{LagF2}
\dot x&=&\sqrt{\epsilon} y\hspace{1.3in}
\dot \beta  = \frac{\epsilon}{2\pi}
\int _0^{2\pi} (\tilde R_\gamma /\sqrt{I \Omega })\sin \phi
d\phi
\nonumber \\
\dot y&=&-\frac{\sqrt{\epsilon}}{2\pi} \int _0^{2\pi}R_0 d\phi _\gamma
\hspace{0.45in}
\dot  I =\frac{-\epsilon}{2\pi} \int _0^{2\pi}(\tilde R_\gamma \sqrt{4I /\Omega }) \cos \phi  d\phi
\end{eqnarray}
where $x=\bar q_0, y=\bar p_0$.
Then, after carrying the averaging (that is similar to \S \ref{Equivalent}),
renaming variables, and setting $\Omega =2\omega$,
we obtain Eq. (\ref{ARE_PR}).

Moreover, the averaged reduced system (without friction) has an effective Hamiltonian
\begin{equation}\label{PRHam}
H_{PR}=\frac{1}{2}y^2+\epsilon \bar M -\epsilon
\frac{\sigma}{2\omega}I+\epsilon  f\frac{I\cos 2\beta}{4\omega}
\end{equation}
that can provide insights on the global phase space structure of this averaged reduced model.
Notice that since angle-action variables are used in our derivation, the Hamiltonian that we have obtained is canonical.  See reference \cite{Ha99} for comparison.

\subsection{Merging Local Bifurcation Analysis with Global Geometry of the Effective Hamiltonian}

Detailed local bifurcation analysis of the averaged reduced system can be used to reveal the ranges of $\sigma $ and other parameters, $f,\mu $, where the  desired dynamics may be available \cite{NaNa93, NaCh95}.  For example, for $f=2.5, \mu =0.5/\omega $ where $\omega =\omega _6=1.17$, the frequency response curves for the averaged reduced system can be depicted as in Figure \ref{Bifurcation}.
These frequency response curves can be obtained by studying
the fixed points of Eq. (\ref{ARE_PR}) and their stability that are parametrized by the detuning parameter $\sigma $.
Notice first that for all fixed points, $y=y_e=0$.  Moreover, there are two types of fixed points, $I=0$ or $I\neq 0$.

\subsubsection{Fixed Points with $I=0$.}  There exist two sub-cases: (i) $x=0$ and the fixed point is $(x,y,q,p)=(0,0,0,0)$, (ii) $x=x_e=\pm 12.59$ and the fixed points are $(x,y,q,p)=(x_e,0,0,0)$.  These two types of solutions can be obtained
by solving the following equation
\begin{equation}\label{MX}
\bar M_x=
\left(\sum _{k=1}^{13}2kc_{2k}(I)x^{2k-2}\right)x=0
\end{equation}
with $I=0$.  Notice that in our study of the fixed points and their stability  with $I=0$, the following Hamiltonian polar coordinates has been used:
\begin{equation}\label{Polar}
q =\sqrt{2I/\omega }\sin \beta ,
\hspace{0.5in}
p =\sqrt{2I\omega }\cos \beta
\end{equation}
That is why $(q,p)=(0,0)$ when $I=0$.

For case (i), the eigenvalues of the fixed point are given by
\begin{equation}
\lambda = \frac{-\epsilon \mu \pm \sqrt{(\epsilon  \mu )^2 -8\epsilon b}}{2},
\hspace{.5in} -\frac{\epsilon \mu }{2} \pm \frac{\epsilon
}{4\omega }
\sqrt{f^2-4(\sigma -2b)^2}
\end{equation}
where $b=c_2(0)=-0.47\times 10^{-4}$.
For example, as $\sigma $ decreases from $3$ to $-2$, this fixed point goes through
two bifurcations, one at $1.1455$ and another at $-1.1457$.  In the interval
$(-1.1457,1.1455)$, it is a rank two saddle (saddle$\times $saddle), marked with red
dashed lines.
In the intervals $(-\infty, -1.1457)$ and $(1.1455, \infty)$, it is a saddle$\times $stable foci, marked with magenta dashed and dotted lines.  See Figure \ref{Bifurcation}(a)(b).

For case (ii), the eigenvalues of the fixed points are given by
\begin{equation}
\lambda = \frac{-\epsilon \mu \pm \sqrt{(\epsilon  \mu )^2 -4\epsilon \bar M_{xx}}}{2},
\hspace{.5in} -\frac{\epsilon \mu }{2} \pm \frac{\epsilon
}{4\omega }
\sqrt{f^2-4(\sigma -2\omega \bar M_I)^2}
\end{equation}
where $\bar M_{xx}$ is the second partial derivative of the averaged Morse potential term $\bar M$ with respect to $x$ and $\bar M_I$ is the first partial derivative of $\bar M$ with respect to $I$.
In the interval $(-2,3)$, these two fixed points $(\pm 12.59,0,0,0)$ are
stable (stable foci$\times$stable foci), marked with blue lines.
See Figure \ref{Bifurcation}(c)(d).
Even though they do go through bifurcations at $55.6726$ and $57.9636$, we do not include their analysis in this paper because they are far outside of the region where we have found the desired dynamics.

\subsubsection{Fixed Points with $I\neq 0$}
\begin{figure}[tb]
\begin{center}
\includegraphics[width=\textwidth]{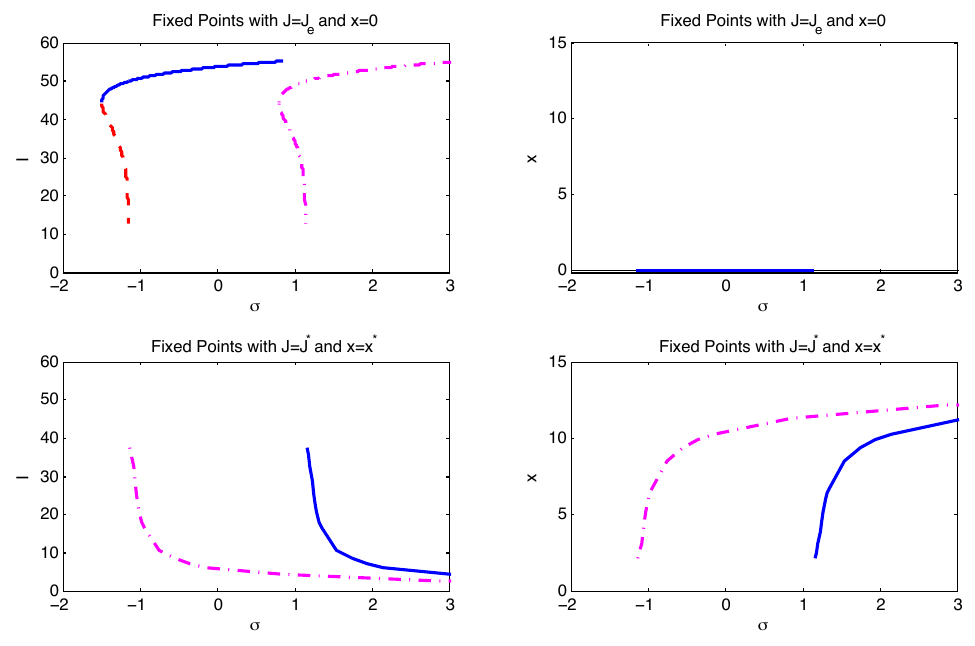}
\end{center}
\caption{\label{Bifurcation2} {\footnotesize
Figures (a)(b) and Figures (c)(d) show the curves of fixed points for
$I_e=I_e(\sigma), x=0$ and $I^*=I^*(\sigma ),x^*=x^*(\sigma )$ respectively.
Figure (a): (1) the fixed points of the upper left branch of $I_e=I_e(\sigma)$
are stable (stable foci$\times$stable foci) for $-1.4825<\sigma <1$, marked with
the blue line ;
(2) the fixed points of lower left branch are saddle$\times$saddle, marked with the red dashed line;
(3) the fixed points of right branch are saddle$\times$stable foci, marked with the
magenta dashed and dotted line.
Figure \ref{Bifurcation2}(b): for simplicity of presentation, only part of the $x=0$ solution that corresponds to upper left branch in Figure (a) has been marked
with the blue line (stable foci$\times$ stable foci).
Figures (c)(d): (1) the fixed points of the left branches
are saddle$\times$stable foci, marked with the magenta dashed and dotted lines; (2) the fixed points of the right branches are stable (stable foci$\times$stable foci), marked with the blue lines.
}}
\end{figure}

For these fixed points,
$\beta=\beta _e=\tfrac{1}{2}\arcsin \sqrt{2\omega \mu /f}$.   Again, there exist two sub-cases.
(iii) x=0 and the fixed points are $(x,y,I,\beta )=(0,0,I_e,\beta _e)$ where
$I_e=I_e(\sigma )$ is obtained by solving
the following equation for $I_e$
\begin{equation}\label{MI}
\bar M_I=\frac{1}{2\omega }\left(\sigma \pm \frac{1}{2}\sqrt{f^2-4(\omega \mu)^2}
\right)
\end{equation}
with $x=0$.  See Figures \ref{Bifurcation2}(a)(b) for the curves of these fixed points.
(iv) $x\neq 0$ and the fixed points are $(x,y,I,\beta )=(x^*,0,I^*,\beta _e)$ where
$I^*=I^*(\sigma ), x^*=x^*(\sigma )$ are obtained by solving the two nonlinear equations
Eq. (\ref{MX}) and Eq. (\ref{MI}) simultaneously.  See Figures \ref{Bifurcation2}(c)(d) for the curves of these fixed points.

For case (iii), the eigenvalues of the fixed points are given by
\begin{equation}
\lambda = \frac{-\epsilon \mu \pm \sqrt{(\epsilon  \mu )^2 -4\epsilon \bar M_{xx}}}{2},
\hspace{.2in} -\frac{\epsilon \mu }{2} \pm \frac{\epsilon
}{2\omega }
\sqrt{(\omega \mu )^2+4\omega fI\bar M_{II} \cos 2\beta _e}.
\end{equation}
Numerical computation shows that (1) the fixed points of the upper left branch of $I_e=I_e(\sigma)$ in Figure \ref{Bifurcation2}(a)
are stable (stable foci$\times$stable foci) for $-1.4825<\sigma <1$;
(2) the fixed points of lower left branch are saddle$\times$saddle;
(3) the fixed points of right branch are saddle$\times$stable foci.
Note: for simplicity of presentation, only part of the $x=0$ solution in Figure \ref{Bifurcation2}(b)
that corresponds to upper left branch in Figure \ref{Bifurcation2}(a) has been marked
with blue line (stable foci$\times$ stable foci).

For case (iv), the characteristic equation of the eigenvalues of the fixed points is given by
\begin{equation}
\lambda ^4 +\tau _1\lambda ^3+\tau _2\lambda ^2+\tau _3\lambda +\tau _4=0
\end{equation}
where
\begin{eqnarray*}
\tau _1&=&
2 \epsilon \mu \\
\tau _2&=&
\epsilon ^2\mu ^2+\epsilon \bar M_{xx}-\epsilon ^2fI\bar M_{II}\cos 2\beta /\omega\\
\tau _3&=&
\epsilon ^2\mu \bar M_{xx} -\epsilon ^3\mu fI\bar M_{II}\cos 2\beta /\omega \\
\tau _4&=&
\epsilon ^3fI(\bar M_{Ix}-\bar M_{xx}\bar M_{II})\cos 2\beta /\omega .
\end{eqnarray*}
Numerical computation shows that (1) the fixed points of the left branches of
Figures \ref{Bifurcation2}(c)(d) are saddle$\times$stable foci; (2) the fixed points of the right branches are stable (stable foci$\times$stable foci).

\begin{figure}[tb]
\begin{center}
\includegraphics[width=\textwidth]{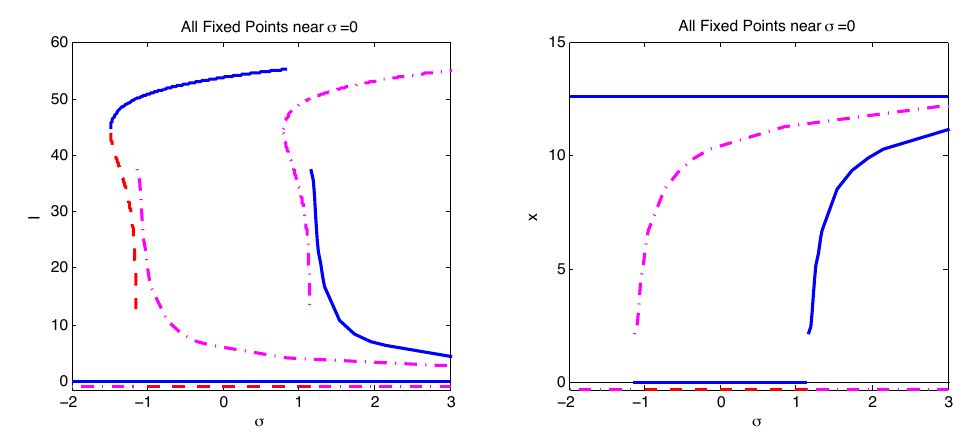}
\end{center}
\caption{\label{Bifurcation} {\footnotesize
Frequency response curves.
Figure shows
how the fixed points and their stability change as the detuning parameter $\sigma $ is varied.   Since for all the fixed points, $y=0$ and $\beta =\beta _e$, only
(a) $I$ vs $\sigma $ and (b) $x$ vs $\sigma $ are plotted.  For example, in the neighborhood of $\sigma =(-0.5, 0.5)$,  the system has four types of  fixed points for each $\sigma $:
(1) the  red dashed curves where $I=0, x=0$ denote saddle $\times $ saddle;
(2) the solid blue curves where $I=0, x=x_e=\pm 12.59$ denote
stable foci$\times$stable foci (only positive $x_e$ is drawn);
(3) the solid blue curves where $I=I_e, x=0$ denote stable foci$\times$stable foci;
(4) the magenta dashed dotted curves where $I=I^*,x=x^*$ denote saddle$\times$ stable foci.
}}
\end{figure}

Figure \ref{Bifurcation} is the combined result of the case studies.
The frequency response curves in this figure show how the fixed points and their stability change as the detuning parameter $\sigma $ (and hence the frequency $\Omega $) of parametric excitation is varied.   By analyzing the relationship between
all these curves, we can see that  the averaged
reduced model may have the desired dynamics in the neighborhood of $\sigma =0$.
For example, at $\sigma =0$, the phase space has four types of fixed points as follows:
\begin{enumerate}
\item
$(x,y,q,p)=(0,0,0,0)$ is a rank two saddle.
\item
$(x,y,q,p)=(x_e,0,0,0)$ where $x_e=\pm 12.59$ are stable foci.  Notice that these fixed points locate at the bottom of the potential well of the averaged reduced system.
\item
$(x,y,I,\beta )=(0,0,I_e,\beta _e)$ where $I_e=45.74\times \omega _6,
\beta _e=\tfrac{1}{2}\arcsin \sqrt{2\omega _6 \mu /f}$  are stable foci.  Notice that the $(x,y)$ coordinates of these fixed points mark the DNA division.
\item
$(x,y,I,\beta )=(x^*,0,I^*,\beta _e)$ where $x^*=\pm 8.7478, I^*=5.312\times \omega _6$ are stable foci $\times$ saddle.  If the frictional coefficient $\mu $ is small.  These fixed points are essentially rank one saddles.
\end{enumerate}
According to the theory of tube dynamics \cite{KoLoMaRo00, GaKoMaRo2005, GaKo2006}, these modified rank one saddles may provide a low energy pathway from the neighborhood of the bottom of the well to
the region which marks the DNA division.

\subsubsection{Global Geometry of the Effective Hamiltonian}

While the local bifurcation analysis does provide many basic ingredients for our study, it does not by itself give a clear and global picture of the dynamics of the averaged reduced system.
Hence, the effective Hamiltonian $H_{PR}(x,y,I,\beta)$ is needed  to fill in this gap.  Notice that in the Hamiltonian polar coordinates Eq. (\ref{Polar}), $\beta =\pi/2$ corresponds to the case where $p=0$.  Therefore,
\begin{equation}\label{Effective}
H_{PR}(x,0,I,\pi/2)=\epsilon \bar M -\epsilon
\frac{\sigma}{2\omega}I-\epsilon  f\frac{I}{4\omega}
\end{equation}
provides an effective potential for the averaged reduced system.

\begin{figure}[tb]
\begin{center}
\includegraphics[width=\textwidth]{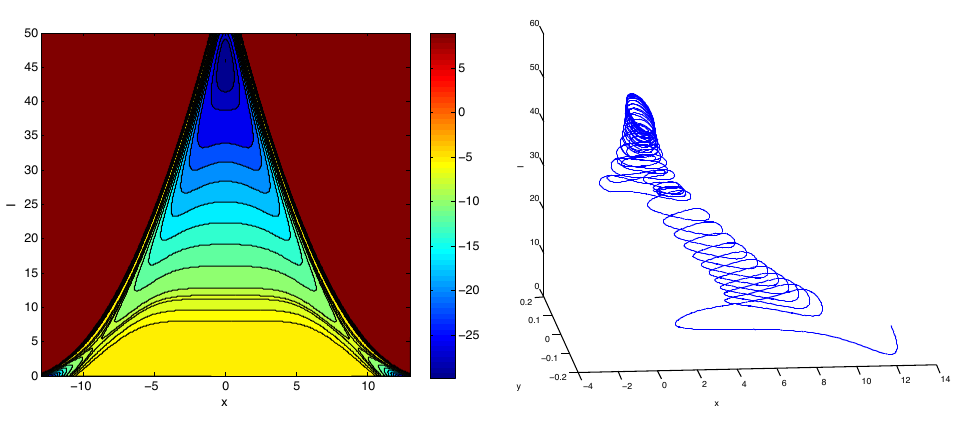}
\caption{\label{Global} {\footnotesize
(a) Contour plots of the effective potential energy in
the $(x,I)$ space.  (b) An example of a trajectory that shows how  the parametric resonance drives the averaged reduced system from its almost equilibrium state to its open state in the $(x,y,I)$ space.
}}
\end{center}
\end{figure}

Figure \ref{Global}(a) shows
the energy contours of this effective potential when $\sigma =0$.  It provides us with the insights for the global phase space structure of the averaged reduced system.  From the figure, we can see clearly how the four types of fixed points
obtained previously fit together within the global geometry of the effective Hamiltonian.  Moreover,
\begin{itemize}
\item
The parametric excitation represented by the parameters $f, \sigma $ has turned $(0,I_e)$, which marks the DNA division, into a sink,
\item
it also creates two
low barrier rank one saddles that are
close to the two DNA equilibrium states $(\pm 12.59,0)$.
\end{itemize}
Hence, the addition of parametric excitation to the averaged reduce system should allow certain trajectories with a little energy in the trigger mode to move from an
almost equilibrium state,
over the saddle, navigate down the energy contours, and reach the sink.  All these insights drawn from the local bifurcation analysis and the effective Hamiltonian  enable us to generate
a class of trajectories that show how the parametric resonance drive the averaged reduced system from its almost equilibrium state to its open state.   Figure \ref{Global}(b) shows an example of this kind of trajectories in the $(x,y,I)$ space.

\subsection{Parametric Resonance Drives DNA to Division}

The data for this trajectory are given as follow.
For the system parameter, we have $f=2.5,\mu=0.5/\omega _6 , \sigma =0$.
For the initial condition, we have $x_0=12.59, y_0=0, I_0=0.54015 \times \omega _6, \beta _0=0$.  The integration is done using the averaged reduced equations,
Eq. (\ref{ARE_PR}).
Notice that this trajectory in Figure \ref{Global}(b) starts at the equilibrium position of the reactive mode but with certain small amount of energy in the trigger mode
($6$th mode).  Without the parametric excitation, the system will liberate near the equilibrium state if there is no friction or die down if the friction exists.  However, if the parametric excitation is turned on at $t=0$ with the data provided above, the $1:2$ parametric resonance will inject energy into the trigger mode, increase the value of $I$, make
the trajectory to reach the region that marks the DNA division ($x=0$) but with large energies in the trigger mode.  See Figure \ref{PendulaC} below for a physical interpretation of this class of trajectories when it is near the DNA open state.

Here, we would like to make a remark on the amount of initial energy in the trigger mode.  Note that
if we increase the amplitude $f$ of the parametric excitation, the magenta dashed dotted curve of Figure \ref{Bifurcation}(a) will shift downward.   Similarly, the rank one saddles in Figure \ref{Global}(a) will also shift downward.  These mean that the amount of initial energy needed in the trigger mode, namely, the value of $I$, can be lower for large $f$.  Numerical simulations of the full system confirm this observation.

\subsection{Extend the Results to the Reduced and the Full Models}

\begin{figure}[tb]
\begin{center}
\includegraphics[width=\textwidth]{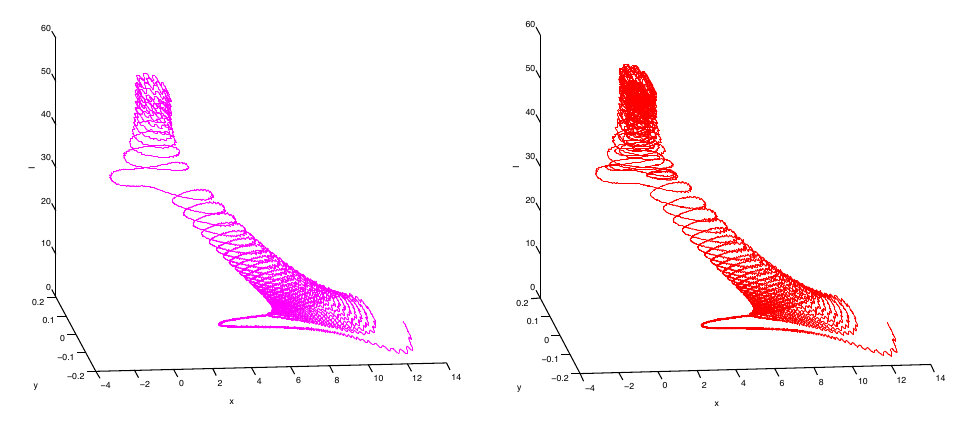}
\caption{\label{Full} {\footnotesize
(a) An example of a trajectory that shows how  the parametric resonance drive the reduced system from its almost equilibrium state to its open state.
(b) A corresponding trajectory for the full model.
}}
\end{center}
\end{figure}

Figure \ref{Full} shows two corresponding trajectories, one for the
reduced model and another for the full model.
The trajectory in Figure \ref{Full}(a) is generated with the following data.
For the system parameters, we set  $f=2.5,\mu=0.5/\omega _6 , \sigma =0$ as before.
For the initial condition, we set $x_0=q_0(0)=12.85,
y_0=p_0(0)=0, I_0=0.8680\times \omega _6, \beta _0=0$ and integrate the trajectory
using the reduced equations, Eq. (\ref{Reduced_PR}) where
\begin{equation}\label{Polar2}
q_6(0)=\sqrt{2I_0/\omega_6 }\sin \beta _0,
\hspace{0.5in}
p_6(0) =\sqrt{2I_0\omega _6 }\cos \beta_0 .
\end{equation}

As for the trajectory in Figure \ref{Full}(b), it  is generated with the following data.
For the system parameters, we set $f=2.5,\mu=0.5/\omega _6, \sigma =0$ as before.
For the initial condition, we set $x_0=q_0(0)=12.85, y_0=p_0(0)=0, I_0=0.86708\times \omega _6, \beta _0=0$ and integrate the trajectory using the full equations, Eq. (\ref{Full_PR}) where $(q_6(0),p_6(0))$ is obtained by Eq. (\ref{Polar2}) and $\theta _k(0)$ are determined by the Fourier modal transformation,
Eq. (\ref{Fourier}).

\begin{figure}[tb]
\begin{center}
\includegraphics[width=\textwidth]{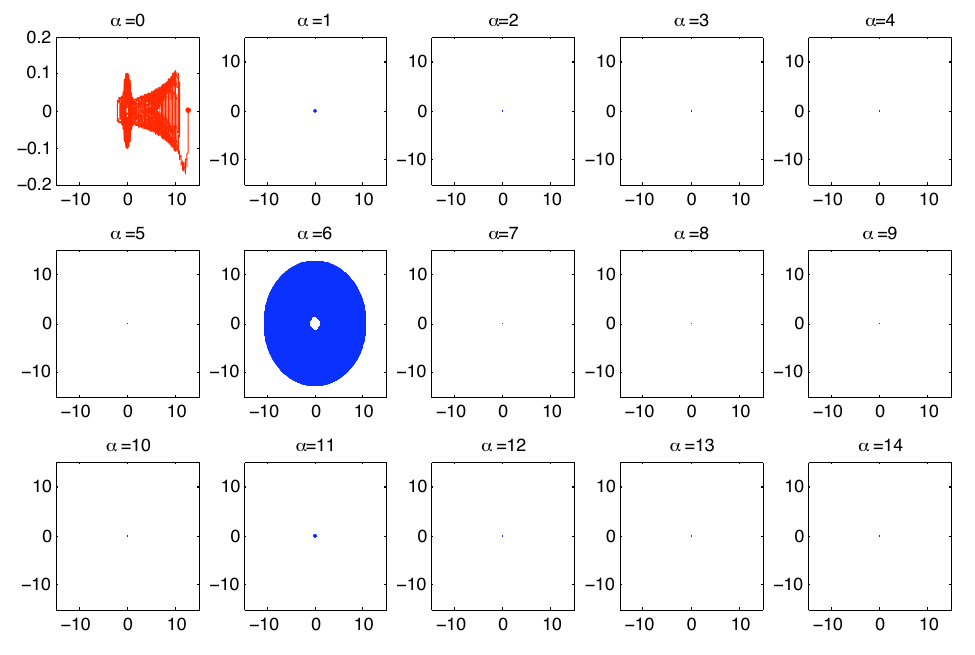}
\end{center}
\caption{\label{15plots_PR} {\footnotesize
Figure shows the projection of the same trajectory of the full model on the phase space of its first 15 modes.  The parametric excitation injects energy into the 6th mode, some of which transfers to the reactive mode and drive the full system from
the almost equilibrium position to the region which marks the DNA division.
}}
\end{figure}

\subsubsection{Remarks on This Class of Special Trajectories}

Here, we would like to make a few remarks:
\begin{itemize}
\item
Despited its simplicity (when compared to the full system), the averaged reduced
model is surprisingly accurate as illustrated by the fact that the three
trajectories in Figure \ref{Global}(b),
Figure \ref{Full}(a), and Figure \ref{Full}(b) are very similar.  Without a careful study of the averaged reduced equations, it may be difficult to guess that such a class of trajectoris will exist in the averaged reduced model, let alone in the reduced and the full model.
\item
It is also interesting to point out that the initial conditions for the reduced and the full models are essentially the same.   Hence, if a small amount of initial energy is in a single Fourier mode,  the dynamics of the reduced model of a two mode truncation
looks  very similar to the dynamics of the full equations of our DNA model.
Figure \ref{15plots_PR} may illustrate this point in another way.   This figure shows the projection of the same trajectory of the full model on the phase space of its first 15 modes.
We observe that (i) the parametric excitation injects energy into the 6th mode, some of which transfers to the reactive mode and drive the full system from the almost equilibrium position to the region that marks the  DNA division, (ii) only an extremely small amount of energy transfers from the excited mode to the other modes.
This observation again shows that the two mode truncation can provide an adequate reduced model for studying the control of DNA division via parametric resonance.
See Figure \ref{15Modes} for comparison.

\begin{figure}[tb]
\begin{center}
\includegraphics[width=\textwidth]{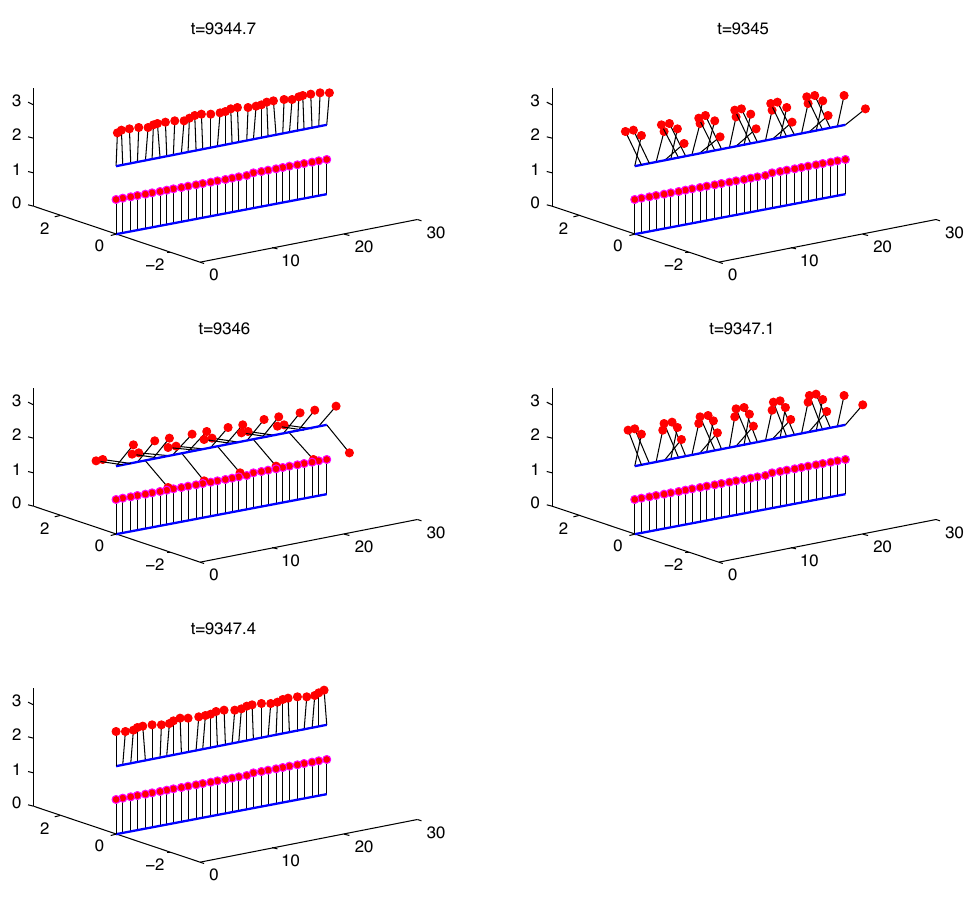}
\end{center}
\caption{\label{PendulaC} {\footnotesize
Figure shows a sequence of 5 snapshots of the evolution of the DNA chain near its opening state
when its corresponding solution trajectory is near the region that marks the DNA division.
}}
\end{figure}

\item
Figure \ref{Full}(b) shows that the trajectory starts at the equilibrium position of the reactive mode but with certain small amount of energy in the trigger mode.  Without the parametric excitation, the system will liberate near the equilibrium state if there is no friction or die down if the friction exists.  However, if the parametric excitation is turned on at $t=0$ with the data provided above, the $1:2$ parametric resonance will inject energy into the trigger mode, increase the value of $I$, make
the trajectory to reach the region that marks the DNA division ($x=0$) but with large energies in the trigger mode.

Figure \ref{PendulaC} shows a sequence of 5 snapshots of the evolution of this solution trajectory in the physical space
when it is near the region that marks the DNA division.
Notice that  the DNA chain that corresponds to this solution is near its
open state ---- the chain is  near the upright position (with
its average angle near zero) but with a periodic swing.
Figures \ref{PendulaC} (a) to (e) show one of
its swings.  Since the 6th mode is used, Figure \ref{PendulaC}(a) shows a curve of pendula
with six peaks (same for the other 4 figures).  Moreover, it is also interesting to point out that the time elapsed
between Figure \ref{PendulaC}(a) to Figure \ref{PendulaC}(e) is 2.7 units of time which is nothing but the
period of the parametric excitation ($2\pi/2\omega$ with $\omega=1.17)$.
\item
Even though qualitatively the trajectories for the averaged reduced system and the reduced system look the same, quantitatively there are certain discrepancy in their initial conditions.  The main reason is that in computing the average of the Morse term $\bar M$, the Taylor expansion at $\theta =0$ is used.  Hence, the equilibrium point for the reactive mode of the averaged reduced system
is given by $(x,y)=(12.59,0)$  instead of $(q_0,p_0)=(\sqrt{n}\theta _e,0)=(12.85,0)$ for  those of the reduced system.  We expect that if we compute $\bar M$ with another Taylor expansion at $\sqrt{n}\theta _e$, the discrepancy will be much less.
See \S \ref{Analytical}.
However, since our concern at this stage is to prove the concept that parametric resonance can be used to control the DNA division, we will not tackle this numerical issue for now.
\item
For the cases where the initial energy is in more than one mode, say in the 6th and the 7th modes, numerical simulation of the full model shows that this kind of trajectories still exist as long as one of the mode is dominant and the parametric excitation is in resonance with the dominant mode.   This should not surprise us because while the parametric resonance will inject energy into the dominant mode, the friction will damp out the other mode.

\item
More studies may be needed in the future for the tradeoffs between
the amplitude $f$, the detuning parameter $\sigma $, the frictional coefficient $\mu $ on one hand and the initial action-phase $I_0,\beta _0$ on the other.
\end{itemize}

\section{Conclusion}

In this paper we have studied  the internal resonance, energy transfer, activation mechanism,
and control of a model of DNA division via parametric resonance.
Our study has been based on a methodology merging geometric reduction, partial averaging, techniques of  chaotic transport,  and control via parametric resonance. This methodology is not limited to our current model and its application can be extended to more general molecular and mechanical systems (possibly involving multiple scales).
This study also highlights the importance of inertial effects in molecular dynamics, such effects are usually ignored in classical studies of molecular systems with Langevin equations (although Langevin equations preserve the Gibbs distribution as an invariant distribution, they do  not account for the electrostatic screening nor the hydrophobic effects of the solvent and they
introduce strong and not necessarily justified assumptions on the dynamic of molecular systems). This is why we have in this first study analyzed a noiseless system. Further studies are required to analyze the effects of  inhomogeneity, helicity (see \cite{Chou84b}), and noise (we note again here that there is no unique way to introduce noise in such systems, henceforth the dynamical aspect of the model may become strongly biased without proper experimental validation).
The inertial effects studied in this paper are important, not only from a general modeling aspect, but also because
 they can be targeted for purpose of control (possibly with low intensity electro-magnetic fields). There is also increasing evidence that these inertial effects play significant biological roles \cite{YaKoMa06}.

It would be  interesting to extend our present framework to the models of DNA with helicity. As has been done in Ref. \cite{Da91,GaRePeDa93}, the effect of the helical geometry of DNA could be incorporated into the present model by introducing additional coupling between every $N$ pendula (nucleotides), where $N$ is typically 4. This kind of coupling could induce another pathway for intramolecular energy transfer, which could in turn make the excitation of the reactive mode even more effective. The helicity of DNA could also be responsible for the coupling between the dynamics of DNA bases and that of the DNA backbone. Since the division of DNA is associated with the slowest-scale dynamics among the dynamics of the bases as we have seen in the present study, the division dynamics of DNA could be coupled (or in resonance) with the slow-scale dynamics of the DNA backbone. In order to study this kind of couplings between the bases and the backbone of DNA, it would be important to implement a novel DNA model that takes into consideration the three-dimensional helical geometry more directly along the lines of Refs. \cite{Yomo83, Vasu09}. Since the helical geometry is ubiquitous among biomolecules, one can expect that the helicity plays a fundamental role in the functions of biomolecules.

\paragraph{acknowledgments}
This work has been supported by the National Science Foundation under award number NSF-CMMI-092600.   We thank Bryan Eisenhower for giving us the permission to use his figure.  W.S.K. thanks Fields Institute for its invitation to the Marsden Conference (July 2012) where he had an opportunity to exchange ideas about this work with other participants.


\appendix

\section{Partial Averaging of Lagrangian Equations is Equivalent to
Partial Averaging of its Hamiltonian for the Reduced Model}\label{Equivalent}

\subsection{Partial Averaging of Lagrangian Equations}

Recall the reduced equations of motion in the Lagrangian form are given by
\begin{eqnarray}\label{0to1b}
\ddot q_0 \hspace{0.5in} & = & -\epsilon M_0(q_0,q_\gamma), \nonumber \\
\ddot q_\alpha + \omega _\alpha ^2 q_\alpha
&=&
-\epsilon M_\gamma (q_0,q_\gamma ).
\end{eqnarray}
They can be rewritten as a first order system as follows:
\begin{eqnarray}
\dot q_0&=&\hspace{0.13in}\sqrt{\epsilon}p_0 \hspace{1.05in} q_\gamma =p_\gamma
\nonumber \\
\dot p_0&=&-\sqrt{\epsilon}M_0 \hspace{1in} p_\gamma =
-\omega ^2_\gamma q_\gamma -\epsilon M_\gamma
\end{eqnarray}
Clearly, the set of equations in the reactive coordinates are already in the Standard Lagrange Form
with a small parameter $\sqrt{\epsilon}$.  Moreover, by using the angle-action variables defined by
\begin{equation}
q_\gamma =\sqrt{2I_\gamma /\omega _\gamma} \sin \phi _\gamma  \hspace{1in}
p_\gamma = \sqrt{2I_\gamma \omega _\gamma} \cos \phi _\gamma
\end{equation}
where $\phi _\gamma  =\omega _\gamma t+\psi _\gamma $, we transform the other set of equations also into the Lagrange Standard Form
\begin{eqnarray}
\dot q_0&=&\hspace{0.13in}\sqrt{\epsilon}p_0 \hspace{1.05in}
\dot \psi _\gamma =\hspace{0.13in}
\epsilon (M_\gamma /\sqrt{2I_\gamma \omega _\gamma })\sin \phi _\gamma
\nonumber \\
\dot p_0&=&-\sqrt{\epsilon}M_0 \hspace{1in}
\dot I_\gamma =-\epsilon M_\gamma \sqrt{2I_\gamma /\omega _\gamma }\cos \phi _\gamma .
\end{eqnarray}
Hence,  we can apply the standard theory of averaging (by averaging $t$ from $0$ to $2\pi/\omega _\gamma $) and obtain the averaged reduced equaitons
\begin{eqnarray}\label{LagF}
\dot x&=&\sqrt{\epsilon} y\hspace{1.3in}
\dot \psi  = \frac{\epsilon}{2\pi}
\int _0^{2\pi} (M_\gamma /\sqrt{2I \omega })\sin \phi _\gamma
d\phi _\gamma
\nonumber \\
\dot y&=&-\frac{\sqrt{\epsilon}}{2\pi} \int _0^{2\pi}M_0 d\phi _\gamma
\hspace{0.45in}
\dot  I =\frac{-\epsilon}{2\pi} \int _0^{2\pi}(M_\gamma \sqrt{2I /\omega }) \cos \phi _\gamma d\phi _\gamma
\end{eqnarray}
where $x=\bar q_0, y=\bar p_0, I=\bar I_\gamma , \psi =\bar \psi _\gamma ,\omega =\omega _\gamma $.

\subsection{Partial Averaging of Its Hamiltonian}

Recall the averaged reduced Hamiltonian is given by
\begin{equation}
\bar H_2=\frac{1}{2} y^2+ \omega I+ \frac{\epsilon }{2\pi}\int _0^{2\pi} M d\phi _\gamma
\end{equation}
where $M(q_0,q_\gamma)$ is a polynomial in $q_0,q_\gamma$.
Its Hamiltonian equation is given by
\begin{eqnarray}\label{HamF}
\dot x&=& y\hspace{1.8in}
\dot \phi =\omega  +  \frac{\epsilon}{2\pi}
\frac{\partial }{\partial I}\int _0^{2\pi} M d\phi _\gamma
\nonumber \\
\dot y&=&-\frac{\epsilon}{2\pi} \frac{\partial }{\partial x} \int _0^{2\pi}M d\phi _\gamma
\hspace{0.58in}
\dot  I =0
\end{eqnarray}

\subsection{Two  Methods are Equivalent for the Reduced Model}
Recall $M_0=\partial M/\partial q_0, M_\gamma =\partial M/\partial q_\gamma$.   Since
\begin{equation}
\frac{\partial M}{\partial \phi _\gamma}=\frac{\partial M}{\partial q_\gamma }
\frac{\partial q_\gamma }{\partial \phi _\gamma}
\hspace{1in}
\frac{\partial M}{\partial I}=\frac{\partial M}{\partial q_\gamma }
\frac{\partial q_\gamma }{\partial I}
\end{equation}
and $\dot \psi =\dot \phi -\omega $,  the $\dot I$ and the $\dot \psi $ equations
of Eq. (\ref{LagF}) and Eq. (\ref{HamF}) are the same.
Moreover, since
\begin{equation}
\int _0^{2\pi}\frac{\partial M}{\partial q_0}d\phi _\gamma =\frac{\partial }{\partial x}\int _0^{2\pi}M d\phi_\gamma .
\end{equation}
the two set of equations for the reactive coordinates are also the same if they are rewritten in the second order forms.

\newpage



\begin{thebibliography}{99}


\bibitem{Ar98}
{\em Dynamical Systems} III, edited by V. Arnold,   Springer-Verlag,
New York, (1998).







\bibitem{BrChKiVe93}
H. W. Broer, S. N. Chow, Y. Kim, and G. Vegter,  A normally elliptic Hamiltonian bifurcation, {\em ZAMP} \textbf{33}, pp.389-432, (1993).



\bibitem{Careri75}
Careri G, Fasella P, Gratton E., Statistical time events in enzymes: a physical assessment, {\em CRC Crit Rev Biochem}, Aug;  3(2), 141-164, (1975).





\bibitem{Chou84}
K. C. Chou, Low-frequency vibrations of DNA molecules,
{\em Biochem J}. 1984 July 1, 221(1), 27-31.

\bibitem{Chou84b}
K. C. Chou, Biological functions of low-frequency vibrations (phonons). III. Helical structures and microenvironment, {\em Biophys J}., 1984 May, 45(5), 881-889.




\bibitem{Da91}
T. Dauxois, Dynamics of breather modes in a nonlinear ``helicoidal'' model
of DNA. \emph{Phys. Lett.} \textbf{A-159}, 390-395 (1991).



\bibitem{Devoe62}
H. DeVoe and I. Tinoco , {\em J. Mol. Biol.}, \textbf{4}, 500 (1962).


\bibitem{duMeMa09}
P. Du Toit, I. Mez‹\'c, and J. Marsden, {\em Physica D} \textbf{238}, 490 (2009).


\bibitem{EiMe07}
B. Eisenhower and I. Mezi\'c, {\em Proceedings of the 46th IEEE
Conference on Decision and Control}, New Orleans 2007, 3976-3981.

\bibitem{EiMe08}
B. Eisenhower, I. Mezi\'c, Actuation requirements in high dimensional oscillator
systems, in {\em Proceedings of the American Control Conference}, 2008.

\bibitem{Ei09}
B. Eisenhower, Targeted escape in large oscillator networks, PhD Dissertation, UCSB, 2009.

\bibitem{EiMe10} B. Eisenhower and I. Mezi\'c, {\em Phys. Rev. E} \textbf{81}, 026603 (2010).

\bibitem{Englander80}
Englander S. W., Calhoun D. B., Englander J. J., Kallenbach N. R., Liem R. K., Malin E. L., Mandal C., Rogero J. R., Individual breathing reactions measured in hemoglobin by hydrogen exchange methods, {\em Biophys J}. 1980 Oct, 32(1), 577-589.





\bibitem{FeLi00}
Z. C. Feng and K. M. Liew, Global bifurcation in parametrically excited system with zero-to-one internal resonance, {\em Nonlinear Dynamics} \textbf{21},
pp. 249-263, (2000).



\bibitem{GaKoMaRo2005}
Gabern, F., W.~S. Koon, J.~E. Marsden, and S.~D. Ross [2005], Theory and
  computation of non-RRKM lifetime distributions and rates in chemical systems
  with three or more degrees of freedom, {\em Physica D} \textbf{211},
  391--406.
\bibitem{GaKo2006}
 Gabern F., W. S. Koon, J. E. Marsden, S. D. Ross, and T. Yanao, [2006],
Application of tube dynamics to non-statistical reaction processes, {\em Few-Body Systems}, Volume \textbf{38}, Numbers 2-4, pages 167 - 172.

\bibitem{GaRePeDa93}
G. Gaeta, C. Reiss, M. Peyrard, and T. Dauxois, Simple models of nonlinear DNA
dyanamics, {\em Rev. Nuovo Cimento}, \textbf{26}, 201-223 (1993)



\bibitem{Ha99}
G. Haller, {\em Chaos Near Resonance}, Springer-Verlag, 1999.



\bibitem{KoLoMaRo00}
W.S. Koon, M.W. Lo, J.E. Marsden, S.D. Ross, Heteroclinic connections between
periodic orbits and resonance transitions in celestial mechanics, {\em Chaos} 10 (2)
(2000) 427-469.


\bibitem{LaZh99}
W. F. Langford and K. Zhan, Interactions of Andronov-Hopf and Bogdanov-Takens bifurcations, {\em Fields Inst. Commun.} \textbf{24}  365-383 (1999)


\bibitem{Letellier1986}
Letellier R., Ghomi M., Taillandier E., Interpretation of DNA vibration modes: I--The guanosine and cytidine residues involved in poly(dG-dC).poly(dG-dC) and d(CG)3.d(CG)3,
{\em J Biomol Struct Dyn}, 3(4), 671--87, (1986)

\bibitem{Lisy96}
Lisy V., Miskovsky P., Schreiber P., On a simple model of low-frequency vibrations in DNA macromolecules, {\em J Biomol Struct Dyn}, 13(4), 707--16, (1996).



\bibitem{LuHoBe96} J. Lumley, P. Holmes, and G. Berkooz, {\em Turbulence, Coherent
Structures, Dynamical Systems and Symmetry}, Cambridge University Press, London, 1996.



\bibitem{Me06}
I. Mezi\'c, On the dynamics of molecular conformation, {\em Proc. Natl. Acad. Sci,} 103
(20) (2006) 7542-7547.


\bibitem{Muto90}
V. Muto, P. S. Lomdahl,  P. L. Christiansen,
Two-dimensional discrete model for DNA dynamics: Longitudinal wave propagation and denaturation, {\em Phys. Rev. A} \textbf{42}  (1990), 7452-7458.



\bibitem{NaNa93}
Nayfeh S. A. and Nayfeh A. H., Nonlinear interactions between two widely spaced modes - external excitation, {\em Int. J. of Bif. Chaos} \textbf{3} (1993) 417-427.

\bibitem{NaNa94}
Nayfeh S. A. and Nayfeh A. H., Energy transfer from high- to low-frequency modes in a flexible structure via modulation, {\em J. Vibr. Acoust.} \textbf{116} (1994) 203-207.

\bibitem{NaCh95}
Nayfeh A. H. and Chin C.-M., Nonlinear interactions in a parametrically  excited system with widely spaced frequencies, {\em Nonlin. Dyn.} \textbf{7} (1995) 195-216.

\bibitem{NaMo95}
Nayfeh, A. H. and Mook, D. T., Energy transfer from high frequency to low frequency modes in structures, {\em Trans. ASME} \textbf{186} (1995) 186-195.





\bibitem{SaVeMu10} J. A. Sanders, F. Verhulst, and J. Murdock, {\em Averaging Methods in Nonlinear Dynamical Systems}, (2010) Springer, 434 pages.

\bibitem{Sobell79}
Sobell H. M., Lozansky E. D., Lessen M., Structural and energetic considerations of wave propagation in DNA, {\em Cold Spring Harb Symp Quant Biol.,} \textbf{43} (Pt 1), 11-19, (1979).

\bibitem{Sobell83}
Sobell H. M., Banerjee A., Lozansky E. D., Zhou G. P., Chou K. C., The role of low-frequency (acoustic) phonons in determining the premelting and melting behaviors of DNA, Structure and Dynamics: {\em Nucleic Acids and Proteins} (Clementi, E. and Sarma, R. H., eds.), pp.181--195, Adenine Press, N.Y.,  (1983).

\bibitem{Song2004}
L. Song, L. Chau, Y. Sakamoto, J. Nakashima, M. Koide, R. S. Tuan,
Electric Field-Induced Molecular Vibration for Noninvasive, High-Efficiency DNA Transfection, {\em Molecular Therapy}, Volume \textbf{9}, issue 4 (April, 2004), 607-616.



\bibitem{TaOwMa11b}
M. Tao, H. Owhadi, and J. E. Marsden,
From Efficient Symplectic Exponentiation of Matrices to Symplectic Integration of High-dimensional Hamiltonian Systems with Slowly Varying Quadratic Stiff Potentials, Appl Math Res Express, 2011-2.

\bibitem{TuVe03}
J M Tuwankotta and F Verhulst  {Hamiltonian system with widely separated frequencies}, {\em Nonlinearity} , \textbf{16} (2003), pp.~689--706.




\bibitem{Vasu09}
V. Vasumathi and M. Daniel. Base-pair opening and bubble transport in a dna double helix
induced by a protein molecule in a viscous medium. Phys. Rev. E, \textbf{80}, 061904 (2009).


\bibitem{Ver00}
F. Verhulst, {\em Nonlinear Differential Equations and Dynamical Systems}, second edition,
Springer, 2000.

\bibitem{Ya04}
L. Yakushevich, {\em Nonlinear Physics of DNA}, Wiley-VCH, 2004.

\bibitem{YaKoMa06}
T. Yanao, W. S. Koon, and J. E. Marsden, {\em Physical Review A} \textbf{73} 052704 (2006).

\bibitem{YaKoMaKe07}
T. Yanao, W.S. Koon, J.E. Marsden, I.G. Kevrekidis, Gyration-radius dynamics
in structural transitions of atomic clusters, {\em J. Chem. Phys.} \textbf{126} (12) (2007)
124102.

\bibitem{YaKoMa09}
T. Yanao, W. S. Koon, and J. E. Marsden, {\em J. Chem. Phys.} \textbf{130}, 144111 (2009).

\bibitem{YaKoMa10}
T. Yanao, W. S. Koon, and J. E. Marsden,
A Nonequilibrium Rate Formula for Collective Motions of Complex Molecular Systems, {\em AIP Conf. Proc.} \textbf{1281}, 1597 (2010).


\bibitem{Yomo83}
S. Yomosa, {\em Phys. Rev. A} \textbf{27}, 2120 (1983).



\bibitem{Zhang89}
C.-T. Zhang, Harmonic and subharmonic resonances of microwave absorption in DNA, {\em Phys. Rev. A} \textbf{40}, 2148-2153 (1989).



\bibitem{ZhBaTu02}
W. Zhang, R. Baskaran, and K. T. Turner, Effect of cubic nonlinearity on auto-parametrically amplified resonant MEMS mass sensor, {\em Sensors and Actuators A} \textbf{102} (2002).

\bibitem{Zhou81}
Zhou G. P.,  Vibrational energy of ringlike DNA molecules,
{\em Shengwu Huaxue Yu Shengwu Wuli Jinzhan}, \textbf{5}, 19--22, (1981).

\bibitem{ZoRa02}
R. S. Zounes and R. H. Rand, Subharmonic resonance in the non-linear Mathieu equation, {\em Int. J. Nonl. Mech.} \textbf{37} (2002).



\end{thebibliography}
\end{document}